\documentclass[12 pt]{article}
\usepackage{amsmath} 
\usepackage{amsfonts,mathrsfs}
\usepackage{amssymb}
\usepackage{color}
\usepackage{graphicx}
\usepackage{epsfig}
\usepackage{subfigure}

{\bf}{}
\newtheorem{lemma}{Lemma}{\bf}{}
\newtheorem{theorem}{Theorem}{\bf}{}
{\bf}{}
\newtheorem{assumption}{Assumption}{\bf}{}
\newenvironment{proof}[1][Proof]{\noindent\textbf{#1.} }{\ \rule{0.5em}{0.5em}}

\def\vmean{v^{\rm mean}}
\def\vmax{v^{\rm max}}

\def\dist{{\rm dist}}

\newcommand{\EXP}[1]{\mathsf{E}\!\left[#1\right] }

\title{Dynamic Coalitional TU Games: Distributed Bargaining among Players' 
Neighbors}

\author{Dario Bauso\footnote{D. Bauso is with Dipartimento di
Ingegneria Informatica, Universit\`a di Palermo, V.le delle
Scienze, 90128 Palermo, ITALY (e-mail: \texttt{dario.bauso@unipa.it})},  Angelia Nedi\'c\footnote{A. Nedi\'c is with the Industrial and Enterprise Systems Engineering Department, University of
Illinois at Urbana-Champaign, Urbana IL 61801 (e-mail: \texttt{angelia@uiuc.edu})
}}

\begin{document}
\maketitle
\begin{abstract}
We consider a sequence of transferable utility (TU) games where, at each time,
the characteristic function is a random vector with realizations restricted to 
some set of values. The game differs from other ones in the literature on 
dynamic, stochastic or interval valued TU games as it combines dynamics 
of the game with an allocation protocol for the players that dynamically 
interact with each other. The protocol is an iterative and decentralized 
algorithm that offers a paradigmatic mathematical description of negotiation 
and bargaining processes. The first part of the paper contributes to 
the definition of a robust (coalitional) TU game and the development of 
a distributed bargaining protocol.
We prove the convergence with probability 1 of the 
bargaining process to a random allocation that lies in the core of the robust
game under some mild conditions on the underlying communication graphs.
The second part of the paper addresses the more general case where 
the robust game may have empty core. In this case, with the dynamic game 
we associate a dynamic average game by averaging over time the sequence of 
characteristic functions. Then, we consider an accordingly  modified 
bargaining protocol. Assuming that the sequence 
of characteristic functions is ergodic and the core of 
the average game has a nonempty relative interior,
we show that the modified bargaining protocol converges with probability 1 to a
random allocation that lies in the core of the average game.
\end{abstract}

\section{Introduction}\label{sec:introduction}
Coalitional games with transferable utilities (TU) have been introduced by 
von Neumann and Morgenstern~\cite{VM}. A coalitional TU game constitutes 
of a set of players, who can form coalitions, and a characteristic function 
that provides a value for each coalition. The value of a coalition  
can be thought of as a monetary value that can be divided 
among the members of the coalition according
to some appropriate fairness allocation rule. 

Coalitional TU games have been used to model cooperation in supply chain 
or inventory management applications~\cite{Drechsel10,HDS00}, as well as in 
network flow applications to coordinate flows of 
goods, materials, or resources between different production/distribution 
sites~\cite{BBP10}. Recently, coalitional games have also sparked 
much interest in communication networks. We refer an interested reader to
tutorial~\cite{SHDHB09} for an in-depth discussion on 
potential applications in this field. 

In this paper, we consider a sequence of coalitional TU games for a finite 
set of players. We assume that the game is played repeatedly over time 
thus generating a sequence of time varying characteristic functions. 
We refer to such a repeated game as \textit{dynamic coalitional TU game}. 
In addition, we assume that a player need not necessarily observe 
the allocations of all the other players, but he can only observe the
allocations of his neighbors, where the neighbors may change in time.

In this setting, we consider robust and average TU games, and for each of 
these games we propose an allocation process to reach a solution of 
the dynamic game. More specifically, we formulate a dynamic bargaining process
as a dynamically changing constrained problem. Iterations are simply bargaining
steps free of any concrete reward allocation. We assume that the reward 
allocation occurs only once at the end of the bargaining process, when 
the players have agreed on how to allocate rewards among themselves. Our main 
objective is to explore distributed agreement on solutions in the core of 
the game, where the players interact only with their neighbors.

In particular, we consider bargaining protocols assuming that each player~$i$ 
obeys rationality and efficiency by deciding on an allocation vector which 
satisfies the value constraints of all the coalitions that include player~$i$. 
This set is termed \textit{bounding set of player~$i$}. At every iteration, 
a player $i$ observes the allocations of some of his neighbors. This is modeled
using a directed graph with the set of players as the vertex set and a 
time-varying edge set composed of directed links $(i,j)$ whenever player~$i$
observes the allocation vector proposed by player~$j$ at time $t$. We refer 
to this directed graph as players' {\it neighbor-graph}. Given a player's 
neighbor-graph, each player~$i$ negotiates allocations by adjusting 
the allocations he received from his neighbors through weight assignments. 
As the balanced allocation may violate his rationality constraints 
(it lies outside player~$i$'s bounding set), the player selects a new allocation
by projecting the balanced allocation on his bounding set. We propose such 
bargaining protocols for solving both the robust and the average TU game.
For each of these games, we use some mild assumptions on the connectivity 
of the players' neighbor-graph and the weights that the players use when 
balancing their own allocations with the neighbors' allocations.
Assuming that the core of the robust game is nonempty, we show that our 
bargaining protocol converges with probability~1 to a common (random) 
allocation in the core. In the case when the core of the robust game 
is empty, we consider an average game that can provide a meaningful solution 
under some conditions on the sequence of the characteristic functions.
Specifically, in this case, we consider a dynamic average game by averaging 
over time the sequence of characteristic functions. We then modify accordingly 
the bounding set and the associated bargaining process by using the dynamic 
average game. This means that the value constraints are defined by using 
the time-averaged sequence rather than the sequence of random characteristic 
functions. Under the assumptions that the time-averaged sequence is ergodic and
that the core of the average game has a nonempty relative interior, we show 
that the players' allocations generated by our bargaining protocol converge 
with probability~1 to a common (random) allocation in the core of the average 
game.

The main contributions of this paper are in the two-fold dynamic aspect of the
games, and in the development and the analysis of distributed 
bargaining process among the players. The novelty in the dynamic aspect
is in the consideration of games with time-varying characteristic functions
and the use of time-varying local players interactions.

The work in this paper is related to stochastic cooperative games 
\cite{G77,SB99,TBT03}. 
However, we deviate from this stochastic framework in a number of aspects 
among which are the existence of a local players interactions captured 
by players' neighbor-graph and the presence of multiple iterations in
the bargaining process, and the consideration of robust game.
Bringing dynamical aspects into the framework of coalitional TU games 
is an element in common with papers~\cite{FP00,H75}. However, 
unlike~\cite{FP00,H75}, the values of the coalitions in this paper are realized
exogenously and no relation is assumed between consecutive samples.
Dynamic robust TU games have been considered in~\cite{BR10} and~\cite{BT09} 
but for a continuous time setting in the former work and for a centralized allocation process in the latter one. 
Convergence of allocation processes is a main topic in~\cite{C98,L02}. 
There, rewards are allocated by a game designer repeatedly in 
a centralized manner and based on the current excess rewards of 
the coalitions (accumulated reward up to current time minus the value of 
the coalition). Our approach, however, differs from that of~\cite{C98,L02} 
as we resort to a \emph{decentralized scheme} where the allocation process is 
the result of a bargaining process among the players with local interactions. 
Convergence of bargaining processes has also been explored under dynamic 
coalition formation~\cite{AS02} for a different dynamic model, where players 
decide both on which coalition to address and what payoff to announce. 
Coalitions form when announced payoffs are all met for all players in the 
coalition. Our bargaining procedure is different as it produces an agreement 
resulting  in the formation of grand coalition, which is stable with respect 
to any sub-coalitions.   

The work in this paper is also related to the literature on agreement among 
multiple agents, where an underlying communication graph for the agents and 
balancing weights have been used with some 
variations~\cite{Tsitsiklis84,NOOT09} to reach an agreement on common 
decision variable, as well as in~\cite{Nedic09,NOP10,Sundhar08c,Sundhar09}
for distributed multi-agent optimization. 

This paper is organized as follows. In Section~\ref{sec:dynamic_game}, we 
introduce the dynamic TU game, the robust game and the bargaining protocol for 
this game. We then motivate the game and give some preliminary results. In 
Section~\ref{sec:main}, we prove the convergence of the bargaining protocol to 
a point in the core of the robust game with probability~1. In 
Section~\ref{sec:average_game}, we introduce the dynamic average game and 
the modified bargaining protocol. We also prove the convergence of this 
protocol to an allocation in the core of the average game with probability~1. 
In Section~\ref{sec:numerical}, we report some numerical simulations to 
illustrate our theoretical study, and we conclude in Section~\ref{conclusions}.

\noindent
{\bf Notation}. 
We view vectors as columns. For a vector $x$, we use $x_i$ or $[x]_i$ to denote
its $i$th coordinate component. 
For two vectors $x$ and $y$, we use $x<y$ ($x\le y$) to denote
$x_i<y_i$ ($x_i\le y_i$) for all coordinate indices $i$. 
We let $x'$ denote the transpose of a vector $x$, and $\|x\|$ 
denote its Euclidean norm.  An $n\times n$ matrix $A$ is row-stochastic if 
the matrix has nonnegative entries $a_{ij}$ and $\sum_{j=1}^n a_{ij}=1$
for all $i=1,\ldots,n$. For a matrix $A$, we use $a_{ij}$ or $[A]_{ij}$
to denote its $ij$th entry. A matrix $A$ is doubly stochastic if
both $A$ and its transpose $A'$ are row-stochastic.
Given two sets $U$ and $S$, we write $U\subset S$ to denote that $U$
is a proper subset of $S$.
We use $|S|$ for the cardinality of a given finite set $S$.

We write $P_X[x]$ to denote the projection of a vector $x$ on a set $X$,
and we write $\dist(x,X)$ for the distance from $x$ to $X$, i.e., 
$P_X[x] = \arg \min_{y \in X} \|x - y\|$ and $\dist(x,X)=\|x-P_X[x]\|$,
respectively.
Given a set $X$ and a scalar $\lambda\in\mathbb{R}$, the set $\lambda X$
is defined by $\lambda X\triangleq\{\lambda x\mid x\in X\}$.
Given two sets $X,Y\subseteq\mathbb{R}^n$, the set sum $X+Y$ is defined by
$X+Y\triangleq\{x+y\mid x\in X,\ y\in Y\}$.
Given a set $N$ of players and a function $\eta: S\mapsto
\mathbb{R}$ defined for each nonempty coalition $S\subseteq N$, 
we write $<N,\eta>$ to denote the transferable utility (TU) game with 
the players' set $N$ and the characteristic function $\eta$.  
We let $\eta_S$ be  the value $\eta(S)$ of the characteristic 
function $\eta$ associated with a nonempty coalition $S\subseteq N$.
Given a TU game $<N,\eta>$, we use $C(\eta)$ to denote the core of the game, 
$C(\eta)=\left\{x \in\mathbb{R}^{|N|} \,\Big|\, \sum_{i\in N}x_i=\eta_N,\  
\sum_{i\in S} x_i\ge\eta_S \hbox{ for all nonempty } S\subset N\right\}.$

\section{Dynamic TU Game and Robust Game}\label{sec:dynamic_game}

In this section, we discuss the main components of a dynamic TU game
and a robust game. We state some basic assumptions on these games 
and introduce a bargaining process conducted by the players to reach an 
agreement on their allocations. 
We also provide some examples of the TU games motivating our 
development and establish some basic results that we use later on in 
the convergence analysis of the bargaining process.  

\subsection{Problem Formulation and Bargaining Process}\label{sec:problem}
Consider a set of players $N=\{1,\ldots,n\}$ and the set of all
(nonempty) \emph{coalitions} $S\subseteq N$ arising among these players.
Let $m=2^n-1$ be the {\it number of possible coalitions}.
We assume that the time is discrete and use $t=0, 1,2,\ldots$ to index the
time slots.

We consider a dynamic TU game, denoted $<N,\{v(t)\}>$, where $\{v(t)\}$ is a
sequence of characteristic functions.
Thus, in the dynamic TU game $<N,\{v(t)\}>$, the players 
are involved in a sequence of instantaneous TU games whereby, at 
each time $t$, the {\it instantaneous TU game} is
$<N,v(t)>$ with $v(t)\in\mathbb{R}^m$ 
for all $t\ge0$. Further, we let $v_S(t)$ denote {\it the value assigned to
a nonempty coalition} $S\subseteq N$ in the instantaneous game $<N,v(t)>$.
In what follows, we deal with dynamic TU games where 
each characteristic function $v(t)$ is a random vector with realizations
restricted to some set of values.
We will consider the robust TU game in this section and Section~\ref{sec:main},
while in Section~\ref{sec:average_game} we deal with the average TU game.

In our development of robust TU game, 
we impose restrictions on each realization of $v(t)$ for every~$t\ge0$.
Specifically, we assume that the grand coalition 
value $v_N(t)$ is deterministic for every~$t\ge0$, while the values $v_S(t)$ 
of the other coalitions $S\subset N$ have a common upper bound. 
These conditions are formally stated in the following assumption.

\begin{assumption}\label{assum:bounded}
There exists $v^{\max}\in\mathbb{R}^m$ such that for all $t\ge0$,
\[v_N(t)=v_N^{\max},\]
\[v_S(t)\le v_S^{\max}\qquad
\hbox{for all nonempty coalitions }S\subset N.\]
\end{assumption}

We refer to the game $<N,v^{\max}>$ as {\it robust game}. 
In the first part of this paper,  we rely on the assumption that the 
robust game has a nonempty core.

\begin{assumption} \label{assum:core}
The core $C(v^{\max})$ of the robust game is not empty.
\end{assumption}

An immediate consequence of Assumptions~\ref{assum:bounded} 
and~\ref{assum:core} is that the core
$C(v(t))$ of the instantaneous game is always not empty. This follows
from the fact that $C(v^{\rm max})\subseteq C(\eta)$ for any $\eta$ satisfying
$\eta_N=v^{\rm max}_N$ and $\eta_S\le v^{\rm max}_S$ for $S\subset N$, 
and the assumption that the core $C(\vmax)$ is not empty.

Throughout the paper, we assume that each player~$i$ 
is {\it rational and efficient}. This translates to each player $i\in N$ 
choosing his allocation vector within the set of allocations satisfying value 
constraints of all the coalitions that include player~$i$. This set is referred
to as the {\it bounding set of player $i$} and, for a generic game $<N,\eta>$, 
it is given by  
\begin{align*}
X_i(\eta)=\left\{x\in\mathbb{R}^n\mid  
\sum_{j\in N}x_j=\eta_N,\ \sum_{j\in S} x_j \geq \eta_S\ 
\mbox{for all $S\subset N$ s.t.~}i\in S\right\}.
\end{align*}
Note that  each $X_i(\eta)$ is polyhedral.

In what follows, we find it suitable to represent 
the bounding sets and the core in alternative equivalent forms. 
For each nonempty coalition $S\subseteq N$, 
let $e_S\in\mathbb{R}^n$ be the incidence vector for the coalition $S$, 
i.e., the vector with the coordinates given by
\[[e_S]_i=\left\{\begin{array}{ll}
1 & \hbox{if $i\in S$,}\cr
0 & \hbox{else}.\end{array}\right.\]
Then, the bounding sets and the core can be represented as follows:
\begin{align}\label{eqn:xipolyh}
X_i(\eta)=\{x\in\mathbb{R}^n\mid e'_N x=\eta_N, \ 
e_S'x\ge \eta_S \ \hbox{for all $S\subset N$
with $i\in S$}\},\end{align}
\begin{align}\label{eqn:corepolyh}
C(\eta)=\{x\in\mathbb{R}^n\mid e'_Nx = \eta_N, \ 
e_S'x\ge \eta_S \ \hbox{for all nonempty $S\subset N$}\}.\end{align}
Furthermore, observe that the core $C(\eta)$ of the game coincides 
with the intersection of the bounding sets $X_i(\eta)$ of all 
players $i\in N=\{1,\ldots,n\}$, i.e.,
\begin{align}\label{eqn:core_eta}
C(\eta) = \cap_{i=1}^n X_i(\eta).\end{align}

We now discuss the bargaining protocol where 
repeatedly over time each player $i\in N$ submits an allocation vector that 
the player would agree on. The allocation vector proposed by player $i$ at 
time $t$ is denoted by $x^i(t)\in \mathbb R^n$, where the $j$th component
$x_j^i(t)$ represents the amount that player $i$ would assign to player $j$. 
To simplify the notation in the dynamic game $<N,\{v(t)\}>$, 
we let $X_i(t)$ denote the bounding set of player $i$ for the 
instantaneous game $<N,v(t)>$, i.e., for all $i\in N$ and $t\ge0$,
\begin{align}\label{eq:bs}X_i(t)
=\left\{x\in\mathbb{R}^n\mid  \sum_{j\in N} x_j = v_N(t),\ 
\sum_{j\in S} x_j \geq v_S(t)
\mbox{ for all $S\subset N$ s.t.~}i \in S\right\}.
\end{align}

\begin{figure} 
\centering
\includegraphics[width=0.34\linewidth]{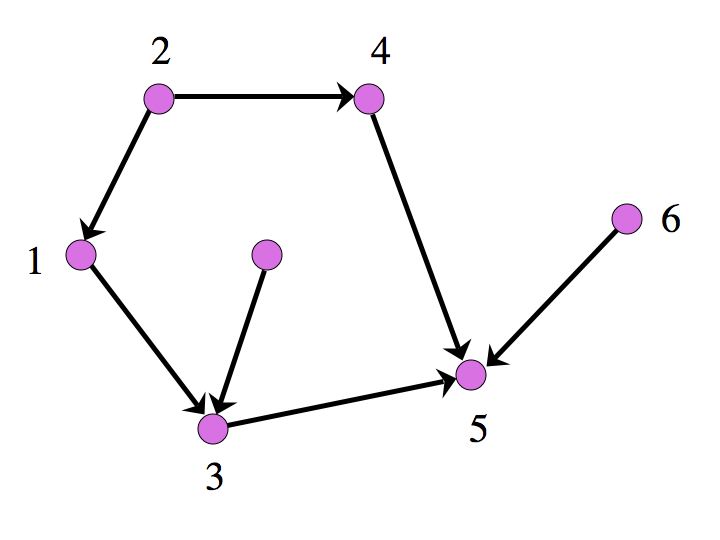}
\qquad\qquad
\includegraphics[width=0.34\linewidth]{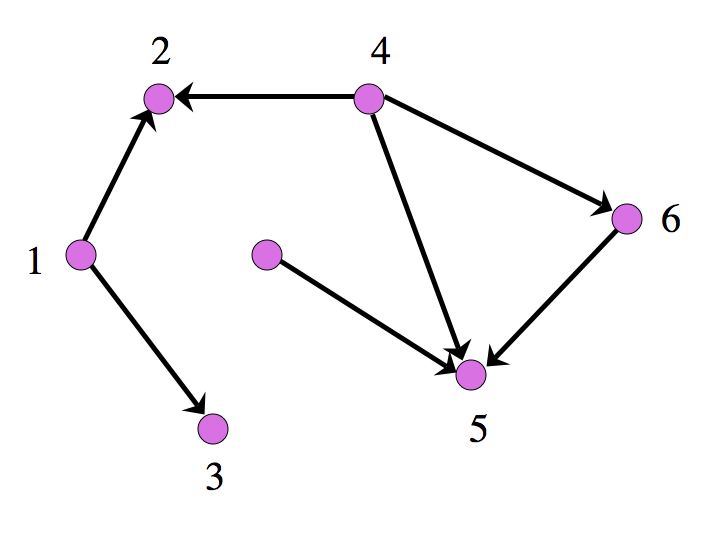}
\caption{Players' neighbor graphs for 6 players
and two different time instances.}\label{fig:graphs}
\end{figure}

We assume that each player may observe the allocations of 
a subset of the other players at any time, which are 
termed as the {\it neighbors of the player}. 
The players and their neighbors at time $t$
can be represented by a directed graph $\mathcal{G}(t)=(N,\mathcal{E}(t))$, 
with the vertex set $N$ and the set $\mathcal{E}(t)$ of directed links.
A link $(i,j)\in \mathcal{E}(t)$ exists if player $j$ is a neighbor 
of player~$i$ at time $t$. 
We always assume that $(i,i)\in\mathcal{E}(t)$ for all $t$, 
which is natural since every player~$i$ can always access its 
own allocation vector. We refer to graph $\mathcal{G}(t)$ as a 
{\it neighbor-graph} at time $t$. In the graph $\mathcal{G}(t)$, 
a player $j$ is a neighbor of player $i$ (i.e., $(i,j)\in\mathcal{E}(t)$) only 
if player $i$ can observe the allocation vector of player~$j$ at time $t$.
Figure~\ref{fig:graphs} illustrates how the players' neighbor-graph 
may look at two time instances.

Given the players' neighbor-graph $\mathcal{G}(t)$,
each player $i$ negotiates allocations by averaging his allocation and 
the allocations he received from his neighbors. 
More precisely, at time $t$, the bargaining process for each player~$i$ 
involves the player's individual bounding set $X_i(t)$, its own 
allocation $x^i(t)$ and the observed allocations $x^j(t)$ of some of 
his neighbors $j$. Formally, we let $\mathcal{N}_i(t)$ be the set of 
neighbors of player $i$ at time $t$ (including himself), i.e.,
\[\mathcal{N}_i(t)=\{j\in N\mid (i,j)\in \mathcal{E}(t)\}.\]
With this notation, the bargaining process is given by:
\begin{equation}\label{eqn:bargaining_neighbor} 
x^i(t+1)= P_{X_i(t)} \left[\sum_{j\in \mathcal{N}_i(t)} a_{ij}(t) x^j(t)
\right]\qquad\hbox{for all $i\in N$ and $t\ge0$},
\end{equation}
where $a_{ij}(t)\ge0$ is a scalar weight that player $i$
assigns to the proposed allocation $x^j(t)$ of player $j$ and 
$P_{X_i(t)}[\cdot]$ is the projection onto 
the player $i$ bounding set $X_i(t)$. The initial allocations $x^i(0),$ 
$i=1,\ldots,n,$ are selected randomly and independently of $\{v(t)\}$.

The bargaining process in~\eqref{eqn:bargaining_neighbor}
can be written more compactly by introducing zero weights 
for players $j$ whose allocations are not available to player $i$ at time $t$.
Specifically by defining $a_{ij}(t)=0$ for all $j\not\in\mathcal{N}_i(t)$
and all $t\ge0$, we have the following equivalent representation of 
the bargaining process:
\begin{equation}\label{eqn:bargaining} 
x^i(t+1) = P_{X_i(t)}\left[\sum_{j=1}^n a_{ij}(t) x^j(t)\right]
\qquad\hbox{for all $i\in N$ and $t\ge0$}.
\end{equation}
Here, 
$a_{ij}(t)=0$ for all $j\not\in\mathcal{N}_i(t)$ while
$a_{ij}(t)\ge0$ for $j\in\mathcal{N}_i(t)$.

We now discuss the specific assumptions on the weights $a_{ij}(t)$ and 
the players' neighbor-graph that we will rely on. 
We let $A(t)$ be the weight matrix with entries $a_{ij}(t)$. 
We will use the following assumption for the weight matrices.

\begin{assumption}\label{assum:weights}
Each matrix $A(t)$ is doubly stochastic with positive diagonal. 
Furthermore, there exists a scalar $\alpha>0$ such that 
\[a_{ij}(t)\ge \alpha\qquad\hbox{whenever}\quad a_{ij}(t)>0.\]
\end{assumption}
In view of the construction of matrices $A(t)$, we see that
$a_{ij}(t)\ge\alpha$ for $j=i$ and perhaps for some players $j$ that
are neighbors of player $i$. The requirement that the positive weights
are uniformly bounded away from zero is imposed to ensure that
the information from each player diffuses to his neighbors in the network 
persistently in time. The requirement on the doubly stochasticity of 
the weights is used to ensure that in a long run each player has equal 
influence on the limiting allocation vector.

It is natural to expect that the connectivity of the players' neighbor-graphs
$\mathcal{G}(t)=(N,\mathcal{E}(t))$ impacts the bargaining process. 
At any time, the instantaneous graph $\mathcal{G}(t)$  
need not be connected. However, for the proper behavior of 
the bargaining process, the union of the graphs 
$\mathcal{G}(t)$ over a period of time is assumed to be connected, 
as given in the following assumption.

\begin{assumption}\label{assum:graph}
There is an integer $Q\ge1$ such that the graph
$\left(N,\bigcup_{\tau=tQ}^{(t+1)Q-1}\mathcal{E}(\tau)\right)$ is 
strongly connected for every $t\ge0$.
\end{assumption}

Assumptions~\ref{assum:weights} and~\ref{assum:graph} together
guarantee that the players communicate sufficiently often to ensure
that the information of each player is persistently diffused over the network 
in time to reach every other player. Under these assumptions, we will study 
the dynamic bargaining process in~\eqref{eqn:bargaining}. We want to provide 
conditions under which the process converges to an allocation in the core of 
the robust game.  
Before this, we provide some motivating examples in the following section.

\subsection{Motivations}\label{sec:motivations}
Dynamic coalitional games capture coordination in a number of network 
flow applications. Network flows model flow of goods, materials, or
other resources between different production/distribution sites~\cite{BBP10}. 
We next provide two examples. The first one is more practically oriented and 
describes a supply chain application \cite{Drechsel10}.
The second one is more theoretical and re-casts 
the problem at hand in terms of controlled flows over hyper-graphs.

\subsubsection{Supply chain}\label{sec:FT}
\begin{figure}[ht]
\centering
\subfigure[Truck leaving ${\bf v}_0$, serving ${\bf v}_1$ and returning
to ${\bf v}_0$.]{
\includegraphics[scale=.44]{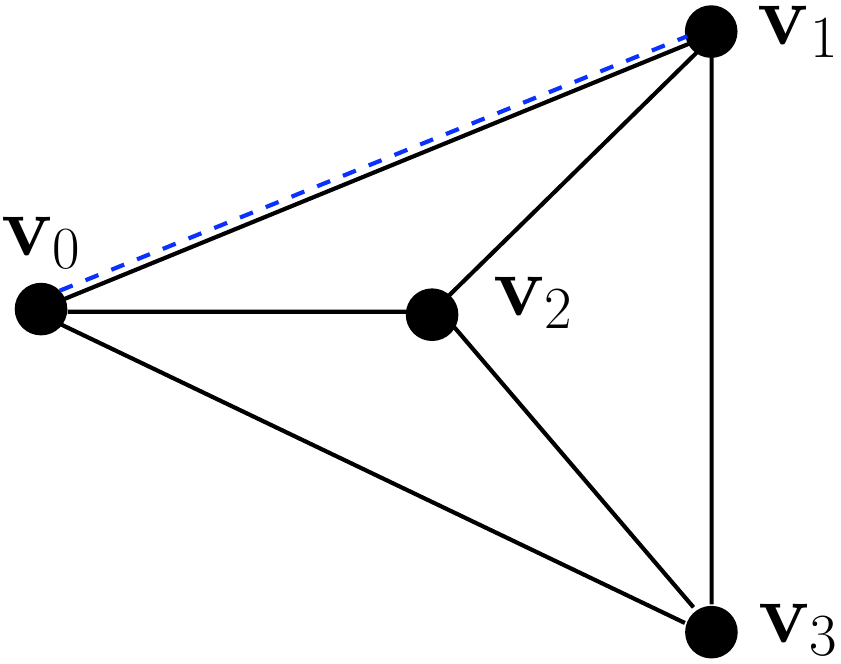}
\label{fig:subfig1}
}
\hskip 1cm
\subfigure[Truck leaving ${\bf v}_0$, serving ${\bf v}_1$ and ${\bf v}_2$, and
then returning to ${\bf v}_0$.]{
\includegraphics[scale=.44]{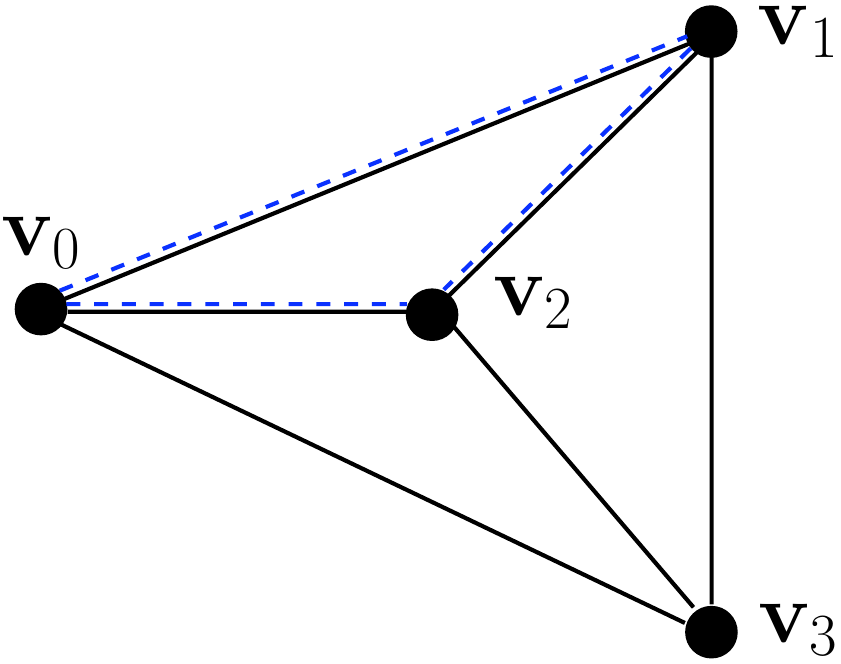}
\label{fig:subfig2}
}
\hskip 1cm
\subfigure[Truck leaving ${\bf v}_0$, serving ${\bf v}_1$,
${\bf v}_2,$ and ${\bf v}_3$, and then 
returning to ${\bf v}_0$.]{
\includegraphics[scale=.44]{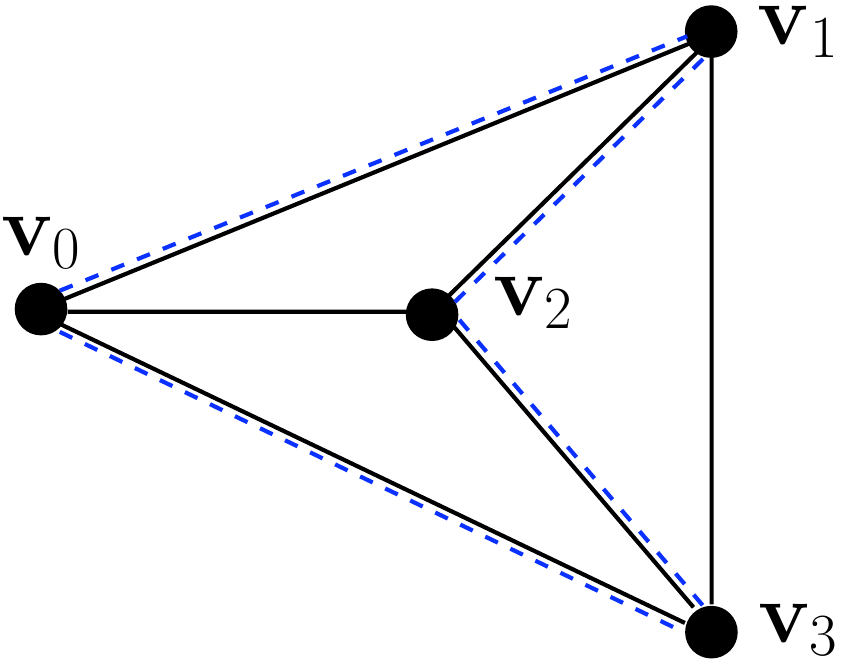}
\label{fig:subfig3}
}
\label{fig:subfigureExample}
\caption[Optional caption for list of figures]{
\small Example of one warehouse 
${\bf v}_0$ and three retailers ${\bf v}_1$, ${\bf v}_2$ and ${\bf v}_3$. 
}
\end{figure}

A single warehouse ${\bf v}_0$ serves a number of retailers 
${\bf v}_i$, $i=1,\ldots,n$, each one facing a demand $d_i(t)$ 
unknown but bounded by pre assigned values  $d_i^{\min} \in \mathbb{R}$ and
$d_i^{\max} \in \mathbb{R}$ at any time period $t\ge0$.    
After demand $d_i(t)$ has been realized, 
retailer ${\bf v}_i$ must choose to either fulfill the demand or not. 
The retailers do not hold any private inventory and, therefore, 
if they wish to fulfill their demands, they must
reorder goods from the central warehouse. Retailers benefit from joint reorders
as they may share the total transportation cost $K$ (this cost could also 
be time and/or players dependent). In particular, if retailer ${\bf v}_i$ 
``plays'' individually, the cost of reordering coincides with the full
transportation cost $K$. Actually, when necessary a single truck will serve 
only him and get back to the warehouse. This is illustrated by the dashed 
path $({\bf v}_0,{\bf v}_1,{\bf v}_0)$ in the network of 
Figure~\ref{fig:subfig1}. The cost of not reordering is 
the cost of the unfulfilled demand $d_i(t)$. 

If two or more retailers ``play'' in a coalition, 
they agree on a joint decision (``everyone reorders'' or
``no one reorders''). The cost of reordering for the coalition
also equals the total transportation cost that must be
shared among the retailers. In this case, when necessary a single truck 
will serve all retailers in the coalition and get back to the warehouse. 
This is illustrated, with reference to coalition $\{{\bf v}_1,{\bf v}_2\}$ 
by the dashed path $({\bf v}_0,{\bf v}_1,{\bf v}_2,{\bf v}_0)$ in
Figure~\ref{fig:subfig2}. A similar comment applies to the coalition 
$\{{\bf v}_1,{\bf v}_2,{\bf v}_3\}$ and the path 
$({\bf v}_0,{\bf v}_1,{\bf v}_2,{\bf v}_3,{\bf v}_0)$ in 
Figure~\ref{fig:subfig3}.
The cost of not reordering is the sum of the unfulfilled 
demands of all retailers. How the players will share 
the cost is a part of the solution generated by the bargaining process.

The cost scheme can be captured by a game with the set 
$N=\{{\bf v}_1,\ldots,{\bf v}_n\}$ of players where 
the cost of a nonempty coalition $S\subseteq N$ is  given by
$$c_S(t)=\min\left\{K,\sum_{i\in S} d_i(t)\right\}.$$
Note that the bounds on the demand $d_i(t)$ reflect into the bounds 
on the cost as follows: for all nonempty $S\subseteq N$ and $t\ge0$, 
\begin{equation}
\label{eq:cs} \min\left\{ K,\sum_{i\in S} d_i^{\min}\right\} 
\leq c_S(t)\leq \min\left\{K,\sum_{i\in S} d_i^{\max}\right\}.\end{equation}
To complete the derivation of the coalitions' values we need to compute the 
cost savings $v_S(t)$ of a coalition $S$ as the
difference between the sum of the costs of the coalitions of the
individual players in $S$ and the cost of the coalition itself, namely,
\[v_S(t) =\sum_{i\in S}c_{\{i\}}(t)-c_S(t).\]

Given the bound for $c_S(t)$ in~\eqref{eq:cs},
the value $v_S(t)$ is also bounded, as given: for any $S\subset N$ and $t\ge0$,
\begin{eqnarray*}
v_S(t) \leq \sum_{i\in S} \min\left\{K, d_i^{\max}\right\}
-\min\left\{K,\sum_{i\in S}d_i^{\min}\right\}.
\end{eqnarray*}
Thus, the cost savings (value) of each coalition is
bounded uniformly by a maximum value. 

\subsubsection{Network controlled flows}\label{sec:FT}
\begin{figure} [htb]
\centering
\includegraphics[width=5cm]{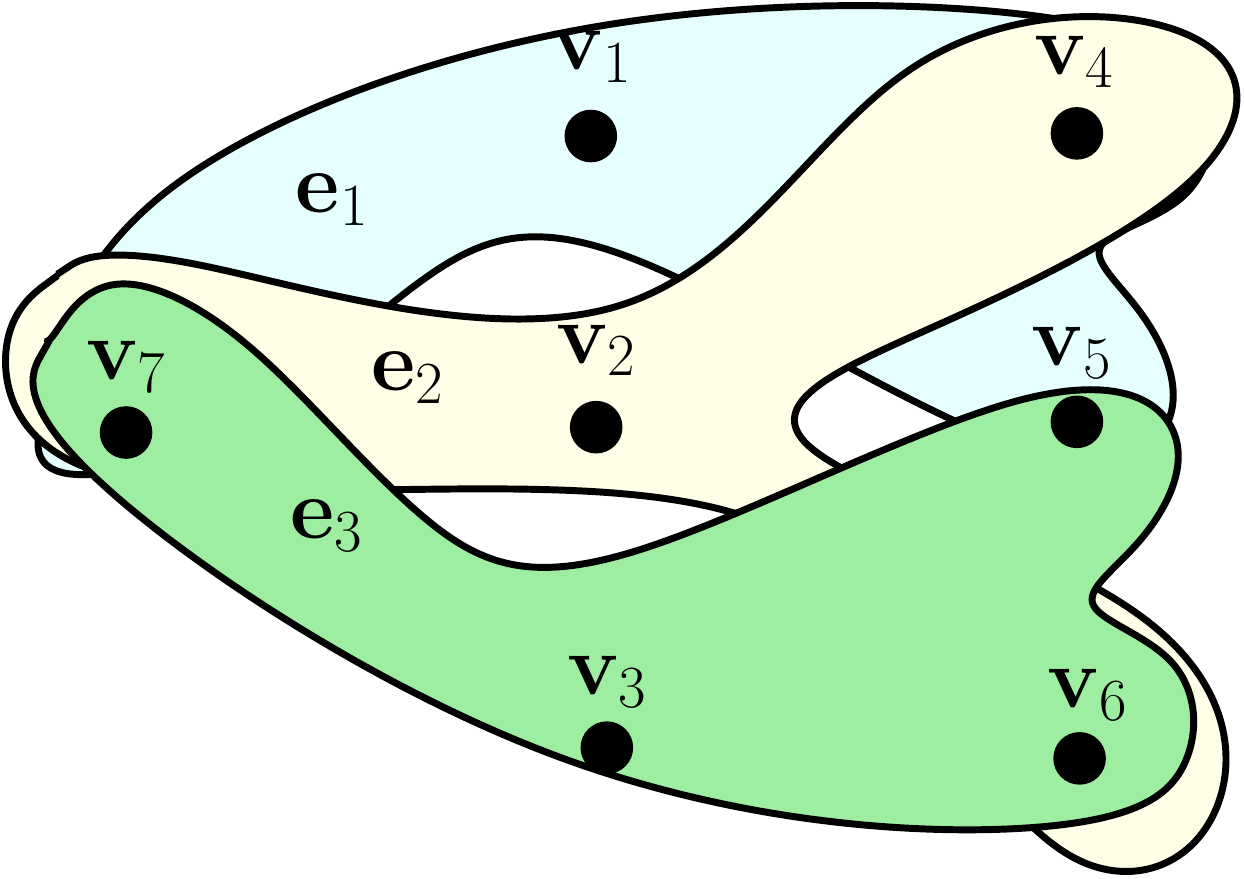}
\caption{Hyper-graph $\mathcal{H}=\{V,E\}$ for a 3-player coalitional game.}  
\label{fig:sat}
\end{figure}

Here, we provide a more theoretical example where the bargaining process
provides a way for the nodes of a hyper-graph to agree on the edge controlled 
flows~\cite{BBP10}.  Consider the hyper-graph $\mathcal{H}=\{V,E\}$ 
with the vertex set $V=\{ {\bf v}_1,\ldots,{\bf v}_m\}$ and
the edge-set $E=\{{\bf e}_1,\ldots,{\bf e}_n\}.$ 
The vertex set $V$ has one vertex per each coalition, 
while the edge-set $E$ has one edge per each player.
An edge ${\bf e}_i$ is incident to a vertex ${\bf v}_j$ if 
the player~$i$ is in the coalition associated with ${\bf v}_j$. 
For a 3-player coalitional game, the vertex-coalition correspondence 
is shown in Table~\ref{Tabb1}. The corresponding 
hyper-graph is shown in Figure~\ref{fig:sat}.

\begin{table}
	\centering
		\begin{tabular}{|c|c|c|c|c|c|c|c|}
\hline		  & ${\bf v_1}$	& ${\bf v_2}$ & ${\bf v_3}$ & ${\bf v_4}$ 
& ${\bf v_5}$ & ${\bf v_6}$ & ${\bf v_7}$\\\hline
	$S$ & $\{1\}$ & $\{2\}$ & $\{3\}$ & $\{1,2\}$ & $\{1,3\}$ & $\{2,3\}$ 
& $\{1,2,3\}$ 	\\\hline
		\end{tabular}\caption{Vertex-coalition correspondence 
for a 3-player coalitional game.}\label{Tabb1}
\end{table}

The incidence relations are described by a matrix 
$B_{\mathcal H}=[e_S']_{S \subseteq N} \in \mathbb R^{m \times n}$ whose 
rows are the characteristic vectors of all coalitions $S\subseteq N$.
The flow control reformulation arises naturally by viewing
an allocation value $\tilde a_i(t)$ of player $i$ 
as the flow on the edge ${\bf e}_i$ associated with player $i$. 
The demand $d_j(t)$ at a vertex 
${\bf v}_j$ is the value $v_S(t)$ of the coalition $S$ associated with the
vertex ${\bf v}_j$. In view of this, an allocation in the core translates into 
satisfying in excess the demand at the vertices, or formally: 
\begin{align*}\label{in} \tilde a(t) \in C(v(t)) \quad 
\Leftrightarrow \quad B_{\mathcal H} \tilde a(t) \geq  d(t).
\end{align*}
In this case, the efficiency condition $\sum_{i\in N} \tilde a_i(t)= v_N(t)$  
corresponds to the equality constraint $\sum_{i\in N} \tilde a_i(t)=d_m(t)$. 
When $d_m(t)$ is a constant, say $d_m(t)=f$ for an $f>0$,
the resulting constraint $\sum_{i\in N} \tilde a_i(t)=f$ can be interpreted as 
a \emph{total flow constraint}, stating that the total flow in the network  
has to be equal to the given value $f$. 

In this framework, the bargaining process is an algorithmic and 
distributed mechanism that ensures the players reach an agreement on 
controlled flows satisfying uncertain demand at the nodes.     

\subsection{Preliminary Results}\label{sec:preliminary}
In this section we derive some preliminary results pertinent to the
core of the robust game and some error bounds for polyhedral sets applicable
to the players' bounding sets $X_i(t)$.
We later use these results to establish the convergence of 
the bargaining process in~\eqref{eqn:bargaining}.

In our analysis we often use the following relation that is valid
for the projection operation on a closed convex set $X\subseteq\mathbb{R}^n$:
for any $w\in\mathbb{R}^n$ and any $x\in X$,
\begin{equation}\label{eqn:projection}
\|P_X[w] - x\|^2\le\|w-x\|^2-\|P_X[w]-w\|^2,
\end{equation} 
This property of the projection operation 
is known as a strictly non-expansive projection property 
(see~\cite{Facchinei2003}, volume II, 12.1.13 Lemma on page 1120). 

We next prove a result that relates the distance $\dist(x,C(\eta))$ between a 
point $x$ and the core $C(\eta)$ with the distances $\dist(x,X_i(\eta))$ 
between $x$ and the bounding sets $X_i(\eta)$. This result will be crucial in 
our later development. The result relies on the polyhedrality of the
bounding sets $X_i(\eta)$ and the core $C(\eta)$, as given 
in~\eqref{eqn:xipolyh} and~\eqref{eqn:corepolyh} respectively, and a special
relation for polyhedral sets. This special relation states that for a nonempty
polyhedral set $\mathcal{P}
=\{x\in\mathbb{R}^n\mid a'_\ell x\le b_\ell, \ \ell=1,\ldots, r\},$
there exists a scalar $c>0$ such that 
\begin{equation}\label{eq:hoffman}
\dist(x,\mathcal{P})\le c \sum_{\ell=1}^r \dist(x,H_\ell)\qquad
\hbox{for all }x\in\mathbb{R}^n,
\end{equation}
where $H_\ell=\{x\in\mathbb{R}^n\mid a_\ell'x\le b_\ell\}$ and the scalar 
$c$ depends on the vectors $a_\ell,\ell=1,\ldots,r$ only.
Relation~\eqref{eq:hoffman} has been established by Hoffman~\cite{Hoffman1952} 
and is known as {\it Hoffman bound}.

\begin{figure} [htb]
\centering
\includegraphics[width=0.4\linewidth]{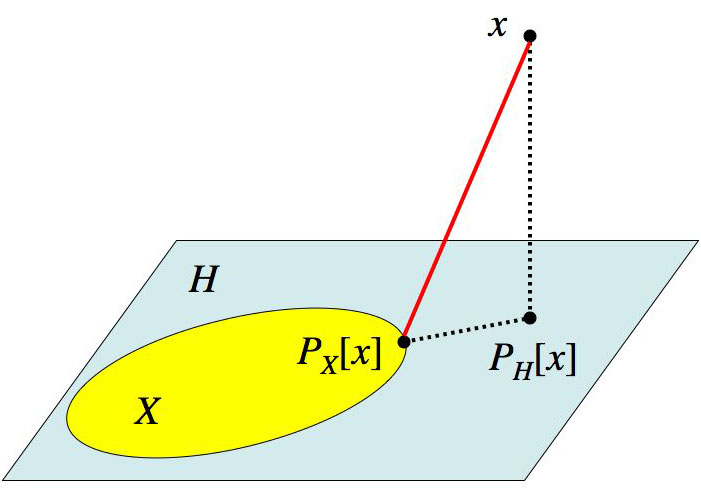}
\caption{Projection on a set $X$ contained in an affine set $H$.}
\label{fig:affine}
\end{figure}

Aside from the Hoffman bound, in establishing the forthcoming 
Lemma~\ref{lemma:polyh}, we also use the fact that the square distance from 
a point $x$ to a closed convex set $X$ contained in an affine set $H$ 
is given by
\begin{equation}\label{eq:affine}
\dist^2(x,X)=\|x-P_H[x]\|^2 + \dist^2(P_H[x], X),
\end{equation}
which is illustrated in Figure~\ref{fig:affine}.

Now, we are ready to present the result relating the values 
$\dist^2(x,C(\eta))$ and $\dist^2(x,X_i(\eta))$. 
\begin{lemma}\label{lemma:polyh}
Let $<N,\eta>$ be a TU game with a nonempty core $C(\eta)$. Then, 
there is a constant $\mu>0$ such that 
\[\dist^2(x,C(\eta))\le \mu\sum_{i=1}^n \dist^2(x,X_i(\eta))
\qquad\hbox{for all $x\in\mathbb{R}^n$},\]
where $\mu$ depends on the collection of vectors 
$\{\tilde e_S\mid S\subset N,\ S\ne\emptyset\}$, where each 
$\tilde e_S$ is the projection of $e_S$ on the hyperplane 
$H =\{x\in\mathbb{R}^n\mid e_N'x=\eta_N\}$.
\end{lemma}

\begin{proof}
Since the hyperplane $H$ contains the core $C(\eta)$ 
(see~\eqref{eqn:corepolyh}), by relation~\eqref{eq:affine} we have
\begin{equation}\label{eq:distance}
\dist^2(x,C(\eta))=\|x-P_H[x]\|^2 + \dist^2(P_H[x], C(\eta))
\qquad\hbox{for all $x\in\mathbb{R}^n$}.
\end{equation}
The point $P_H[x]$ and the core $C(\eta)$ lie in the $n-1$-dimensional affine 
set $H$. By applying the Hoffman bound relative to the affine set $H$ 
(cf.~Eq.~\eqref{eq:hoffman}), we obtain
\[\dist(P_H[x],C(\eta))\le c\, \sum_{S\subset N} \dist(P_H[x], H\cap H_S),\]
where the summation is over nonempty subsets $S$, 
$H_S=\{x\in\mathbb{R}^n\mid e_S'x\ge \eta_S\}$, 
while the constant $c$ depends on the collection $\{\tilde e_S,S\subset N\}$
of projections of vectors $e_S$ on the hyperplane $H$ for $S\subset N$. 
Thus, it follows
\begin{eqnarray}\label{eqn:hofmanninH}
\dist^2(P_H[x],C(\eta))
& \le & c^2\, \left(\sum_{S\subset N} \dist(P_H[x], H\cap H_S)\right)^2\cr
& \le  &c^2 (m-1)\, 
\sum_{S\subset N} \dist^2(P_H[x], H\cap H_S),
\end{eqnarray}
where $m$ is the number of nonempty subsets of $N$ and the last inequality
follows by $(\sum_{j=1}^\ell a_j)^2\le \ell \sum_{j=1}^\ell a_j^2$, 
which is valid for any finite collection of scalars $a_j,j=1,\ldots,\ell$ 
with $\ell\ge1$. Combining Eq.~\eqref{eqn:hofmanninH} with 
equality~\eqref{eq:distance}, we obtain for any $x\in\mathbb{R}^n$,
\begin{eqnarray*}
\dist^2(x,C(\eta))
&\le&  \|x-P_H[x]\|^2 + c^2 (m-1)\, 
\sum_{S\subset N} \dist^2(P_H[x], H\cap H_S)\cr
&\le &  c_1\, \sum_{S\subset N} \left( \|x-P_H[x]\|^2 
+ \dist^2(P_H[x], H\cap H_S)\right),
\end{eqnarray*}
where $c_1=\max\{1,\ c^2 (m-1)\}$. 
Since the set $H$ is affine, in view of relation~\eqref{eq:affine} we have 
$\|x-P_H[x]\|^2 + \dist^2(P_H[x], H\cap H_S)=\dist^2(x,H\cap H_S)$, implying 
that for any $x\in\mathbb{R}^n$,
\[\dist^2(x, C(\eta))\le c_1\, \sum_{S\subset N} \dist^2(x, H\cap H_S).\]

From the preceding relation, it follows for any $x\in\mathbb{R}^n,$
\begin{align}\label{eqn:jedanr}
\dist^2(x, C(\eta))
\le c_1 \sum_{S\subset N} |S|\, \dist^2(x, H\cap H_S),
\end{align}
where $|S|$ denotes the cardinality of the coalition $S$.
Note that 
\begin{align}\label{eqn:dvar}
\sum_{S\subset N} |S|\, \dist^2(x, H\cap H_S)
&=\sum_{S\subset N} \sum_{i\in S}\, \dist^2(x, H\cap H_S)\cr
&=\sum_{i=1}^n \left(\sum_{ \{S\subset N\mid i\in S\} } 
\dist^2(x, H\cap H_S)\right).\end{align}
We also note that $X_i(\eta)\subset H\cap H_S$ 
for each nonempty $S\subset N$ and $i \in S$, which follows by  the definition
of $H_S$ and relation~\eqref{eqn:xipolyh}.
Furthermore, since 
$\dist(x,Y)\le \dist(x,X)$ for any $x\in\mathbb{R}^n$ and 
for any two closed convex sets $X,Y\subseteq\mathbb{R}^n$ such that 
$X\subset Y$, it follows that for all $x\in\mathbb{R}^n$,
\begin{align}\label{eqn:trir}
\dist(x,H\cap H_S)\le \dist(x,X_i(\eta)).\end{align}
By combining relations~\eqref{eqn:jedanr}--\eqref{eqn:trir}, we obtain
\begin{eqnarray*}
\dist^2(x, C(\eta))
& \le & c_1\, \sum_{i=1}^n \left(
\sum_{ \{S\subset N\mid i\in S\} } \dist^2(x,X_i(\eta))\right)\cr
& = &c_1\kappa\, \sum_{i=1}^n  \dist^2(x,X_i(\eta)),
\end{eqnarray*}
where $\kappa$ is the number of coalitions $S$ that contain player $i$, which
is the same number for every player ($\kappa$ does not depend on $i$). 
The desired relation follows by letting $\mu=c_1\kappa$, and
by recalling that $c_1=\max\{1,c^2(m-1)\}$ and that $c$ 
depends on the projections $\tilde e_S$ of vectors $e_S,S\subset N$, 
on the hyperplane $H$.  
\end{proof} 

Note that the scalar $\mu$ in Lemma~\ref{lemma:polyh} does not depend 
on the coalitions' values $\eta_S$ for $S\ne N$. It depends only on the vectors
$e_S$, $S\subseteq N$, and the grand coalition value $\eta_N$.

As a direct consequence of Lemma~\ref{lemma:polyh}, we have the following
result for the instantaneous game $<N,v(t)>$ under the assumptions
of Section~\ref{sec:problem}.

\begin{lemma}\label{lemma:alleta}
Let Assumptions~\ref{assum:bounded} and~\ref{assum:core} hold.
We then have for all $t\ge0$, 
\[\dist^2(x, C(v(t)))
\le \mu \sum_{i=1}^n  \dist^2(x,X_i(t))
\qquad\hbox{for all }x\in\mathbb{R}^n,\]
where $C(v(t))$ is the core of the game $<N,v(t)>$,
$X_i(t)$ is the bounding set of player $i$,
and $\mu$ is the constant from Lemma~\ref{lemma:polyh}.
\end{lemma}
\begin{proof} By Assumption~\ref{assum:core}, we have that 
the core $C(v^{\max})$ is nonempty. Furthermore,
under Assumption~\ref{assum:bounded}, we have $C(v^{\max})\subseteq C(v(t))$ 
for all $t\ge0$, implying that the core $C(v(t))$ is nonempty for all $t\ge0$.

Under Assumption~\ref{assum:bounded}, each core $C(v(t))$ is
defined by the same affine equality corresponding to the grand coalition value,
$e_N'x=v_N^{\max}$.
Moreover, each core $C(v(t))$ is defined through the set of
hyperplanes $H_S(t)=\{x\in\mathbb{R}^n\mid e_S' x\ge v_S(t)\}$, 
$S\subset N$, which have time invariant normal vectors $e_S$, $S\subseteq N$.
Thus, the result follows from Lemma~\ref{lemma:polyh}.
\end{proof}

\section{Convergence to Core of Robust Game}\label{sec:main}
In this section, we prove convergence of the bargaining process
in~\eqref{eqn:bargaining} to a random allocation that lies 
in the core of the robust game with probability~1. 
We find it convenient to re-write the bargaining 
protocol~\eqref{eqn:bargaining} by isolating 
a linear and a non-linear term. The linear term is the vector $w^i(t)$ 
defined as:
\begin{equation}\label{eqn:wit}
w^i(t)=\sum_{j=1}^n a_{ij}(t) x^j(t)\qquad\hbox{for all $i\in N$ and $t\ge0$}.
\end{equation}
Note that $w^i(t)$ is linear in players' allocations $x^j(t)$.
The non-linear term is the error 
\begin{equation}\label{eqn:error1} 
e^i(t)= P_{X_i(t)}[w^i(t)] - w^i(t). 
\end{equation}
Now, using \eqref{eqn:wit} and~\eqref{eqn:error1}, 
we can rewrite protocol~\eqref{eqn:bargaining} as follows:
\begin{equation}\label{eqn:bargaininge} 
x^i(t+1)= w^i(t) + e^i(t)
\qquad\hbox{for all $i\in N$ and $t\ge0$}.
\end{equation}
Recall that the weights $a_{ij}(t)$ are nonnegative and such that
$a_{ij}(t)=0$ for all $j\not\in\mathcal{N}_i(t)$. Also, recall that
$A(t)$ is the matrix with entries $a_{ij}(t)$,
which governs the construction of the vectors $w^i(t)$ in~\eqref{eqn:wit}.

The main result of this section shows that, with probability~1, 
the bargaining protocol~\eqref{eqn:wit}--\eqref{eqn:bargaininge}
converges to the core $C(v^{\max})$ of the robust game $<N,v^{\rm max}>$, 
provided that 
$v(t)=v^{\max}$ happens infinitely often in time with probability~1.
To establish this we use some auxiliary results, which we develop in the next
two lemmas.

The following lemma provides a result on the sequences $x^i(t)$ and shows that 
the errors $e^i(t)$ are diminishing.
\begin{lemma}\label{lemma:basic}
Let Assumptions~\ref{assum:bounded} and~\ref{assum:core} hold.
Also, assume that each matrix $A(t)$ is doubly stochastic.
Then, for the bargaining protocol~\eqref{eqn:wit}--\eqref{eqn:bargaininge}, 
we have 
\begin{itemize}
\item[(a)] The sequence $\left\{\sum_{i=1}^n \|x^i(t+1) - x \|^2\right\}$ 
converges for every $x\in C(v^{\max})$.
\item[(b)] The errors $e^i(t)$ in~\eqref{eqn:error1} are such that 
$\sum_{t=0}^\infty\sum_{j=1}^n \|e^i(t)\|^2<\infty$. In particular,
$\lim_{t\to\infty} \|e^i(t)\|=0$ for all $i\in N$.
\end{itemize}
\end{lemma}

\begin{proof}
By $x^i(t+1)=P_{X_i(t)}[w^i(t)]$ and 
by strictly non-expansive property of the Euclidean projection on 
a closed convex set $X_i(t)$ (see~\eqref{eqn:projection}), 
we have for any $i\in N$, $t\ge0$ and $x\in X_i(t)$,
\begin{align}\label{eqn:reljedan}
\|x^i(t+1) - x \|^2 \leq \|w^i(t) - x \|^2 - \|e^i(t)\|^2.\end{align}
Under Assumptions~\ref{assum:bounded} and~\ref{assum:core},
the core $C(v^{\rm max})$ is contained in the core $C(v(t))$ for all $t\ge0$, 
implying that $C(v^{\max})\subseteq C(v(t))$ for all $t\ge0$. Furthermore,
since $C(v(t))=\cap_{i=1}^n X_i(t)$, it follows that 
$C(v^{\max})\subseteq X_i(t)$ for all $i\in N$ and $t\ge0$.
Therefore, relation~\eqref{eqn:reljedan} holds for all $x\in C(v^{\max})$. 
Thus, by summing the relations in~\eqref{eqn:reljedan} over $i\in N$, we obtain
for all $t\ge0$ and $x\in C(v^{\max})$,
\begin{equation}\label{eq:proj}
\sum_{i=1}^n \|x^i(t+1) - x \|^2 \leq 
\sum_{i=1}^n \|w^i(t) - x \|^2 - \sum_{i=1}^n \|e^i(t)\|^2.
\end{equation}
By the definition of $w^i(t)$ in \eqref{eqn:wit}, using the stochasticity 
of $A(t)$ and the convexity of the squared norm, we obtain
$$\sum_{i=1}^n \|w^i(t) - x\|^2
=\sum_{i=1}^n \left\| \sum_{j=1}^n a_{ij}(t) x^j(t) - x\right\|^2
\leq \sum_{j=1}^n \left(\sum_{i=1}^n  a_{ij}(t) \right) 
\|  x^j(t) - x\|^2.$$
By the doubly stochasticity of $A(t)$, we have
$\sum_{i=1}^n a_{ij}(t)=1$ for every $j$, 
implying $\sum_{i=1}^n \|w^i(t) - x\|^2 
\leq \sum_{i=1}^n \|x^i(t) - x\|^2$. 
By substituting this relation in \eqref{eq:proj}, we arrive at
\begin{equation}\label{eq:ineq}
\sum_{i=1}^n \|x^i(t+1) - x \|^2 
\leq \sum_{i=1}^n \|x^i(t) - x \|^2 - \sum_{i=1}^n \|e^i(t) \|^2.
\end{equation}
The preceding relation shows that the scalar sequence 
$\{\sum_{i=1}^n \|x^i(t+1) - x \|^2\}$ is non-increasing for any 
given $x\in C(v^{\max})$. Therefore, the sequence must be convergent
since it is nonnegative. Moreover, 
by summing the relations in \eqref{eq:ineq} over $t=0,\ldots,s$ and 
taking the limit as $s\rightarrow \infty$, we obtain
$$\sum_{t=0}^{\infty} \sum_{i=1}^n \|e^i(t)\|^2 \leq 
\sum_{i=1}^n \|x^i(0) - x\|^2,$$ which implies that 
$\lim_{t\to\infty}e^i(t)=0$ for all $i\in N$. 
\end{proof}

In our next result, we will use the instantaneous 
average of players allocations, defined as follows:
\[y(t)=\frac{1}{n} \sum_{j=1}^n x^j(t)\qquad\hbox{for all }t\ge0.\]
The result shows that the difference between the bargaining payoff vector 
$x^i(t)$ for any player $i$ and the average $y(t)$ of these payoffs 
converges to 0 as time goes to infinity. The proof essentially uses the line of
analysis that has been employed in~\cite{NOP10}, where the sets $X_i(t)$ 
are static in time, i.e., $X_i(t)=X_i$ for all $t$. In addition, we also use
the rate result for doubly stochastic matrices that has been 
established in~\cite{NOOT09}.

\begin{lemma}\label{lemma:xityt}
Let Assumptions~\ref{assum:weights} and~\ref{assum:graph} hold. Suppose that 
for the bargaining protocol~\eqref{eqn:wit}--\eqref{eqn:bargaininge} we have
\[\lim_{t\to\infty} \|e^i(t)\|=0\qquad\hbox{for all }i\in N.\]
Then, for every player $i\in N$ we have 
\[\lim_{t\to\infty}\|x^i(t)-y(t)\|=0,\qquad
\lim_{t\to\infty}\|w^i(t)-y(t)\|=0.\]
\end{lemma}

\begin{proof}
For any $t\ge s\ge0$, define matrices
\[\Phi(t,s)=A(t)A(t-1)\cdots A(s+1)A(s),\]
with $\Phi(t,t)=A(t)$. Using the matrices $\Phi(t,s)$ and the expression for
$x^i(t)$ in~\eqref{eqn:bargaininge}, we can relate the vectors $x^i(t)$ with 
the vectors $x^i(s)$ at a time $s$ for $0\le s\le t-1$, as follows:
\begin{equation}\label{eqn:short}
x^i(t) = \sum_{j=1}^n [\Phi(t-1,s)]_{ij} \, x^j(s) +\sum_{r=s+1}^{t-1}
\left(\sum_{j=1}^n [\Phi(t-1,r)]_{ij}\, e^j(r-1) \right) +e^i(t-1).
\end{equation}
Under the doubly stochasticity of the matrices $A(t)$, using 
$y(t)=\frac{1}{n}\sum_{j=1}^n x^j(t)$ and relation~\eqref{eqn:short}, we obtain
\begin{equation}
y(t)= {1\over n} \sum_{j=1}^n x^j(s) + {1\over n}\sum_{r=s+1}^{t}
\left(\sum_{j=1}^n  e^j(r-1)\right)\qquad\hbox{for all $t\ge s\ge0$}.
\label{eqn:approxaver}
\end{equation}
By our assumption, we have
$\lim_{t\to\infty} \|e^i(t)\|=0$ for all $i$.
Thus, for any $\epsilon>0$, there is an integer
$\hat s\ge0$ such that $\|e^i(t)\|\le \epsilon$ for all
$t\ge \hat s$ and all $i$.
Using relations~\eqref{eqn:short} and~\eqref{eqn:approxaver} 
with $s=\hat s$, we obtain for all $i$ and $t\ge \hat s+1$,
\begin{eqnarray*}
&&\|x^i(t) -  y(t)\|
= \left\|\sum_{j=1}^n \Big([\Phi(t-1,\hat s)]_{ij} -
{1\over n}\Big)x^j(\hat s) \right.\\
& + & \left. \sum_{r=\hat s+1}^{t-1} \sum_{j=1}^n \Big(
[\Phi(t-1,r)]_{ij} - {1\over n}\Big)e^j(r-1)
+\Big(e^i(t-1)- {1\over n} \sum_{j=1}^n e^j(t-1)\Big)\right\|\\
&\le& \sum_{j=1}^n \Big|[\Phi(t-1,\hat s)]_{ij} - {1\over n}\Big|\
\|x^j(\hat s)\| \\ & + & \sum_{r=\hat s+1}^{t-1} \sum_{j=1}^n
\Big|[\Phi(t-1,r)]_{ij} - {1\over n}\Big|
\|e^j(r-1)\| +  \|e^i(t-1)\| + {1\over n}\sum_{j=1}^n \|e^j(t-1)\|.\\
\end{eqnarray*}
Since $\|e^i(t)\|\le \epsilon$ for all $t\ge \hat s$ and all $i$, 
it follows that
\[\|x^i(t) -  y(t)\|
\le \sum_{j=1}^n \Big|[\Phi(t-1,\hat s)]_{ij} - {1\over n}\Big|\
\|x^j(\hat s)\| +  \epsilon\sum_{r=\hat s+1}^{t-1} \sum_{j=1}^n
\Big|[\Phi(t-1,r)]_{ij} - {1\over n}\Big| \ +  2\epsilon.
\]

Under Assumptions~\ref{assum:weights} and \ref{assum:graph}, 
the following result holds for the matrices $\Phi(t,s)$, 
as shown in \cite{NOOT08} (see there Corollary 1):
\[\left| [\Phi(t,s)]_{ij}-\frac{1}{n}\right|
\le \left(1- \frac{\alpha}{4n^2}\right)^{\left\lceil
\frac{t-s+1}{Q}\right\rceil-2}\qquad\hbox{for all }t\ge s\ge0.\]
Substituting the preceding estimate in the estimate for 
$\|x^i(t) - y(t)\|$, we obtain
\begin{eqnarray*}
\|x^i(t) - y(t)\| 
\le \left(1- \frac{\alpha}{4n^2}\right)^{\left\lceil
\frac{t-\hat s}{Q}\right\rceil-2} \sum_{j=1}^n \|x^j(\hat s)\|  
+ n\epsilon \sum_{r=\hat s+1}^{t-1} 
\left(1- \frac{\alpha}{4n^2}\right)^{\left\lceil
\frac{t-r}{Q}\right\rceil-2}\, +2\epsilon.
\end{eqnarray*}
Letting $t\to\infty$, we see that
\[\limsup_{t\to\infty} \|x^i(t)-y(t)\|
\le n \epsilon\sum_{r=\hat s+1}^{\infty} 
\left(1- \frac{\alpha}{4n^2}\right)^{\left\lceil
\frac{t-r}{Q}\right\rceil-2} +2\epsilon.\]
Note that $\sum_{r=\hat s+1}^{\infty} 
\left(1- \frac{\alpha}{4n^2}\right)^{\left\lceil
\frac{t-r}{Q}\right\rceil-2} <\infty$, 
which by the arbitrary choice of $\epsilon$ yields
\[\lim_{t\to \infty} \|x^i(t)-y(t)\|=0\qquad\hbox{for all $i\in N$}.\]

Now, we focus on $\sum_{i=1}^n\|w^i(t)-y(t)\|$. Since 
$w^i(t)=\sum_{j=1}^n a_{ij}(t)x^j(t)$ and since $A(t)$ is stochastic, 
it follows
$$\sum_{i=1}^n\|w^i(t)-y(t)\|\le
\sum_{i=1}^n\sum_{j=1}^n a_{ij}(t)\|x^j(t) -y(t)\|.$$ 
Exchanging the order of the summations over, and then using 
the doubly stochasticity of $A(t)$, we have
\[\sum_{i=1}^n\|w^i(t)-y(t)\|\le
\sum_{j=1}^n\left(\sum_{i=1}^n a_{ij}(t)\right)\|x^j(t) -y(t)\|
=\sum_{j=1}^n \|x^j(t) - y(t)\|.\]
Since $\lim_{t\to \infty} \|x^j(t) - y(t)\|=0$ for all $j$, it follows that
$\sum_{i=1}^n\|w^i(t) - y(t)\|\to 0,$
thus implying $\|w^i(t) - y(t)\|\to0$ for all $i\in N$.
\end{proof}

Note that Lemma~\ref{lemma:xityt} captures the effects of the matrices 
$A(t)$ that represent players' neighbor-graphs. At the same time,
Lemma~\ref{lemma:basic} is basically a consequence of the projection property 
only. So far, the polyhedrality of the sets $X_i(t)$ has not been used at all. 
We now put all pieces together, namely Lemma~\ref{lemma:alleta} that exploits 
the polyhedrality of the bounding sets $X_i(t)$, Lemma~\ref{lemma:basic} and 
Lemma~\ref{lemma:xityt}. This brings us to the following convergence result
for the robust game $<N,\vmax>$.

\begin{theorem}\label{thm1}
Let Assumptions~\ref{assum:bounded}--\ref{assum:graph} hold.
Also, assume that  
$${\rm Prob}\,\left\{v(t)=v^{\max} \ \ i.o.\right\}=1,$$ 
where $i.o.$ stands for infinitely often.
Then, the players allocations $x^i(t)$ generated by 
bargaining protocol~\eqref{eqn:wit}--\eqref{eqn:bargaininge} 
converge with probability~1 to an allocation in the core  $C(v^{\max})$, i.e., 
there is a random vector $\tilde x\in C(v^{\max})$ 
such that with probability~1,
$$\lim_{t\rightarrow \infty} \|x^i(t)-\tilde x\|=0 \qquad 
\mbox{for all $i\in N$}.$$
\end{theorem}

\begin{proof} 
By Lemma~\ref{lemma:basic}, for each player $i\in N$, the sequence 
$\{\sum_{i=1}^n \|x^i(t)-x\|^2\}$ is convergent for 
every $x\in C(v^{\max})$ and the errors $e^i(t)$ are diminishing,
i.e., $\|e^i(t)\|\to0$. Then, by Lemma~\ref{lemma:xityt} we have
$\|x^i(t)-y(t)\|\to0$ for every $i$. Hence,
\begin{equation}\label{eq:ytpoints}
\hbox{$\{\|y(t)-x\|\}$ is convergent for 
every $x\in C(v^{\max})$}.\end{equation}
We want to show that $\{y(t)\}$ is convergent and that its limit is 
in the core $C(v^{\max})$ with probability~1. For this, we note that
since $x^i(t+1) \in X_i(t)$,  it holds 
\[\sum_{i=1}^n \dist^2\left(y(t+1),X_i(t)\right) 
\leq \sum_{i=1}^n \|y(t+1) - x^i(t+1)\|^2\qquad\hbox{for all $t\ge0$}.\]
The preceding relation and $\|x^i(t)- y(t)\|\rightarrow 0$ for all $i\in N$ 
(cf.~Lemma~\ref{lemma:xityt}) imply 
$$\lim_{t\rightarrow \infty} 
\sum_{i=1}^n \dist^2\left(y(t+1),X_i(t)\right)=0.$$
Under Assumptions~\ref{assum:bounded} and~\ref{assum:core},
by Lemma~\ref{lemma:alleta} we obtain
\[\dist^2\left(y(t+1),C(v(t))\right)\le
\mu \sum_{i=1}^n \dist^2\left(y(t+1),X_i(t)\right)\qquad\hbox{for all $t\ge0$}.
\]
By combining the preceding two relations we see that
\begin{equation}\label{eq:ylimit}
\lim_{t\to\infty}\, \dist^2\left(y(t+1),C(v(t))\right)=0.
\end{equation}

By our assumption, the event 
$\{v(t)=v^{\max}\  \hbox{infinitely often}\}$ happens with probability~1.
We now fix a realization $\{v_\omega(t)\}$ of the sequence $\{v(t)\}$ such that
$v_\omega(t) = v^{\max}$ holds infinitely often (for infinitely many $t$'s).
Let $\{t_k\}$ be a sequence such that 
\[v_\omega(t_k) = v^{\max}\qquad\hbox{for all }k.\]
All the variables corresponding to the realization
$\{v_\omega(t)\}$ are denoted by a  subscript $\omega$.
By relation~\eqref{eq:ytpoints} the sequence $\{y_\omega(t)\}$ is bounded,
therefore $\{y_\omega(t_k)\}$ is bounded. 
Without loss of generality 
(by passing to a subsequence of $\{t_k\}$ if necessary), we assume that
$\{y_\omega(t_k)\}$ converges to some vector $\tilde y_\omega$, i.e., 
\[\lim_{k\to\infty} y_\omega(t_k) = \tilde y_\omega.\] 
Thus, the preceding two relations and Eq.~\eqref{eq:ylimit} imply that 
$\tilde y_\omega\in C(v^{\max})$.
Then, by relation~\eqref{eq:ytpoints}, we have that 
$\{\|y_\omega(t)-\tilde y_\omega\|\}$ is convergent, from which we conclude 
that $\tilde y_\omega$ must be the unique accumulation point of the sequence
$\{y_\omega(t)\}$, i.e., 
\[\lim_{t\to\infty} y_\omega(t)=\tilde y_\omega,\qquad 
\tilde y_\omega\in C(v^{\max}).\]
This and the assumption ${\rm Prob}\,\left\{v(t)=v^{\max} \ \ i.o.\right\}=1,$
imply that the sequence $\{y(t)\}$ converges with probability~1 to a 
random point $\tilde y\in C(v^{\max})$. Since by Lemma~\ref{lemma:xityt} 
we have $\|x^i(t)-y(t)\|\to0$ for every $i$, it follows that the sequences 
$\{x^i(t)\},i=1,\ldots,n,$ converge with probability~1 to a common random 
point in the core $C(v^{\max})$.\end{proof}

\section{Dynamic Average Game}\label{sec:average_game}
When the core of the robust game is empty, the core of the instantaneous 
average game can provide a meaningful solution under some conditions on 
the distribution of the functions $v(t)$. In what follows, we focus on 
the instantaneous average game associated with the dynamic TU game 
$<N,\{v(t)\}>$.  In the next sections, we define the instantaneous 
average game, we introduce a bargaining protocol for the game and investigate 
the convergence properties of the bargaining protocol. 

\subsection{Average Game and Bargaining Protocol}\label{sec:ave_game}
Consider a dynamic TU  game $<N,\{v(t)\}>$ with each $v(t)$ being a random 
characteristic function. With the dynamic game we associate a 
{\it dynamic average game} $<N,\{\bar v(t)\}>$, where $\bar v(t)$ is 
the average of the characteristic functions $v(0),\ldots,v(t)$, i.e.,
\[\bar v(t)=\frac{1}{t+1}\,\sum_{k=0}^{t} v(k)\qquad\hbox{for all }t\ge0.\]
An {\it instantaneous average game at time $t\ge0$} is the game 
$<N,\bar v(t)>$. We let $C(\bar v(t))$ denote the core of the instantaneous 
average game at time~$t$ and let
$\bar X_i(t)$ denote the bounding set of player $i$ for the 
instantaneous game $<N,\bar v(t)>$, i.e., for all $i\in N$ and $t\ge0$,
\begin{equation}\label{eqn:bs_ave}
\bar X_i(t)=\left\{x\in\mathbb{R}^n\mid  e_N'x = \bar v_N(t),\ 
e_S'x \geq \bar v_S(t)
\mbox{ for all $S\subset N$ such that $i \in S$}\right\}.
\end{equation}
Note that $\bar X_i(0)=X_i(0)$ for all $i\in N$ since $\bar v(0)= v(0)$.
In what follows, we {\it assume that $\bar X_i(t)$ are nonempty} 
for all $i\in N$ and all $t\ge0$. 

In this setting, the bargaining process for the players is given by
\begin{equation}\label{eqn:bargaining_ave} 
x^i(t+1) = P_{\bar X_i(t)}\left[\sum_{j=1}^n a_{ij}(t) x^j(t)\right]
\qquad\hbox{for all $i\in N$ and $t\ge0$},
\end{equation}
where $a_{ij}(t)\ge0$ is a scalar weight that player $i$
assigns to the proposed allocation $x^j(t)$ received from player $j$ at time 
$t$. The initial allocations $x^i(0)$, $i\in N$,  are selected 
randomly and independently of $\{v(t)\}$. Regarding the weights $a_{ij}(t)$, 
recall that these weights are reflective of the players' neighbor-graph:
$a_{ij}(t)=0$ for all $j\not\in\mathcal{N}_i(t)$, where
$\mathcal{N}_i(t)$ is the set of neighbors of player $i$ (including himself) 
at time $t$, while we may have $a_{ij}(t)\ge0$ only for $j\in\mathcal{N}_i(t)$.

\subsection{Assumptions and Preliminaries}\label{sec:asum_prelim}
In this section, we provide our assumptions for the average game, and
discuss some auxiliary results that we need later on in the convergence 
analysis of the bargaining protocol. Regarding the random characteristic 
functions $v(t)$ we use the following assumption.

\begin{assumption}\label{assum:ergodic} 
The sequence $\{v(t)\}$ is ergodic, i.e., with probability~1, we have
\[\lim_{t\to\infty}\bar v(t)=v^{\rm mean}\qquad
\hbox{with $v^{\rm mean}\in\mathbb{R}^m$}.\]
Furthermore, $v_N(t)=\vmean$ for all $t$ with probability~1.
\end{assumption}
Assumption~\ref{assum:ergodic} basically says that the grand coalition
value $v_N(t)$ is constant with probability~1.
Note that Assumption~\ref{assum:ergodic} is satisfied, for example, when
$\{v_S(t)\}$ is an independent identically distributed sequence
with a finite expectation  $\EXP{v_S(t)}$ for all $S\subset N$.

We refer to the TU game $<N,v^{\rm mean}>$ as {\it average game}, which
is well defined under Assumption~\ref{assum:ergodic}. We let 
$C(v^{\rm mean})$ be the core of the average game and $\bar X_i$ be 
the bounding set for player $i$ in the game, i.e.,
\[
C(\vmean)=\left\{x\in\mathbb{R}^n\mid e_N'x=v^{\rm mean}_N,\ 
e_S' x\ge v^{\rm mean}_S\quad \hbox{for all nonempty $S\subset N$}
\right\},\]
\begin{equation}\label{eqn:bs_mean}
\bar X_i=\left\{x\in\mathbb{R}^n\mid e_N'x=v^{\rm mean}_N,\ 
e_S' x\ge v^{\rm mean}_S\quad \hbox{for all $S\subset N$ such that $i\in S$}
\right\}.\end{equation}
Note that the average core $C(v^{\rm mean})$ lies in the hyperplane 
$H=\{x\in \mathbb{R}^n\mid e_N'x=\vmean_N\}$. Hence, the dimension of the
core is at most $n-1$. We will in fact assume that the dimension of the average
core $C(v^{\rm mean})$ is $n-1$, by requiring the existence of a point $z$ in 
the core such that all other inequalities defining the core are satisfied 
as strict inequalities. 

Specifically, we make use of the following assumption. 
\begin{assumption}\label{assum:cores}
There exists a vector $\hat z\in C(v^{\rm mean})$ such that
\[e_S'\hat z>\vmean_S\qquad\hbox{for all nonempty }S\subset N.\]
\end{assumption}
Assumption~\ref{assum:cores} basically says that $\hat z$ is in the relative 
interior of the core $C(\vmean)$ of the average game and that
the core $C(\vmean)$ has dimension $n-1$. This assumption will be important 
in establishing the convergence of the bargaining 
protocol~\eqref{eqn:bargaining_ave}. In particular,
the following result will be important, which is an immediate consequence of 
Assumption~\ref{assum:cores} and the polyhedrality of the cores $C(\bar v(t))$.

\begin{lemma}\label{lemma:relint}
Let Assumptions~\ref{assum:ergodic} and~\ref{assum:cores} hold. Then,
with probability 1, for every $z$ in relative interior of $C(\vmean)$ 
there exists $t_z$ large enough such that $z$ is in the relative interior of
$C(\bar v(t))$ for all $t\ge t_z$ with probability~1. 
\end{lemma}

\begin{proof}
Let $z$ be in the relative interior of $C(\vmean)$ which exists by 
Assumption~\ref{assum:cores}. Thus, $e_N'z=\vmean_N$ and 
$e_S'z>\vmean_S$ for all $S\subset N$.
By Assumption~\ref{assum:ergodic}, with probability~1 
we have $\bar v_N(t)=\vmean_N$ and $\bar v_S(t)\to \vmean_S$ for $S\subset N$.
Hence, $\bar v_N(t)'z=(\vmean_N)'z$ for all $t$ with probability~1. 
Furthermore, there exists a random time $t_z$ large enough so that 
with probability 1,
\[e_S'z>\bar v_S(t)\qquad\hbox{for all $S\subset N$ and all $t\ge t_z$},\]
implying that $z$ is in the relative interior of $C(\bar v(t))$ for all
$t\ge t_z$ with probability~1.
\end{proof}
Lemma~\ref{lemma:relint} shows that the sets $C(\bar v(t))$ and 
$C(\vmean)$ have the same dimension for large enough $t$ with probability~1.
In particular,  this lemma implies that the cores $C(\bar v(t))$ are nonempty 
with probability~1 for all $t$ sufficiently large. 

Aside from Lemma~\ref{lemma:relint}, 
in our convergence analysis of the bargaining protocol, we make use of two 
additional well-known results. One of them is the super-martingale convergence 
theorem due to Robbins and Siegmund~\cite{Robbins1971} (it can
also be found in~\cite{Polyak87}, Chapter~2.2, Lemma~11).

\begin{theorem}\label{thm:supermartingale}
Let $\{V_k\}, \{g_k\},$ and $\{h_k\}$ be non-negative random
scalar sequences. Let $F_k$ be the $\sigma$-algebra generated by $V_1,
\ldots, V_k, g_1, \ldots g_k,h_1, \ldots h_k.$ Suppose that almost surely,
\[\EXP{V_{k+1} \mid F_k } \leq V_{k} - g_{k} + h_{k}\qquad \hbox{for all
$k$,} \]
and $\sum_{k} h_k <\infty$  almost surely.  Then, almost surely 
both the sequence $\{V_{k}\}$ converges to a non-negative random variable and  
$\sum_{k} g_{k} < \infty$.
\end{theorem}

The other well-known result that we use is pertinent 
to two nonempty polyhedral sets whose description 
differs only in the right-hand side vector. The result can be found 
in~\cite{Facchinei2003}, 3.2.5~Corollary, pages 258--259.

\begin{lemma}\label{lemma:twopolyhsets}
Let $P_b$ and $P_{\tilde b}$ be two polyhedral sets given by
\[P_b=\{x\in\mathbb{R}^n\mid Bx\le b\},\qquad
P_{\tilde b}=\{x\in\mathbb{R}^n\mid Bx\le \tilde b\},\]
where $B$ is an $m\times n$ matrix and $b,\tilde b\in\mathbb{R}^m$. 
Then, there is a scalar $L>0$ such that for every $b,\tilde b\in\mathbb{R}^m$
for which $P_b\ne\emptyset$ and $P_{\tilde b}\ne\emptyset$, we have
\[\dist(x,P_{\tilde b})\le L \|b-\tilde b\|
\qquad\hbox{for any $x\in P_b$,}\]
where the constant $L$ depends on the matrix $B$.
\end{lemma}

For two polyhedral sets $P_b$ and $P_{\tilde b}$ as in 
Lemma~\ref{lemma:twopolyhsets} we also have 
\begin{equation}\label{eqn:dist}
\dist(x,P_{b})\le \dist(x,P_{\tilde b})+L\|b-\tilde b\|
\qquad\hbox{for any }x\in\mathbb{R}^n
\end{equation}
(see Rockafellar~\cite{Rockafellar1998}, 4.34 Lemma on page 132).

\subsection{Convergence to Core of Average Game}\label{sec:ave_game}
Here, we show the convergence of the bargaining protocol to the core
of the average game. In our analysis, we find it convenient to re-write 
the bargaining protocol~\eqref{eqn:bargaining_ave}
in an equivalent form by separating a linear and a non-linear term. 
The linear term is the vector $\bar w^i(t)$ given by
\begin{equation}\label{eqn:wit_ave}
\bar w^i(t)=\sum_{j=1}^n a_{ij}(t) x^j(t)\qquad\hbox{for all 
$i\in N$ and $t$}.
\end{equation}
The non-linear term is expressed by the error 
\begin{equation}\label{eqn:error2} 
\bar e^i(t)= P_{\bar X_i(t)}[\bar w^i(t)] - \bar w^i(t). 
\end{equation}
Now, using relations~\eqref{eqn:wit_ave} and~\eqref{eqn:error2}, 
we can rewrite~\eqref{eqn:bargaining_ave} as follows:
\begin{equation}\label{eqn:bargaininge_ave} 
x^i(t+1)= \bar w^i(t) + \bar e^i(t)
\qquad\hbox{for all $i\in N$ and all $t\ge0$}.
\end{equation}

We first show some basic properties of the players' allocations
by using the preceding equivalent description of the bargaining 
protocol~\eqref{eqn:bargaining_ave}. These properties hold under
the doubly stochasticity of the weights $a_{ij}(t)$ that comprise the matrix 
$A(t)$.
\begin{lemma}\label{lemma:basic_ave}
Let Assumptions~\ref{assum:ergodic} and~\ref{assum:cores} hold.
Also, let the matrices $A(t)$ be doubly stochastic. Then, for bargaining 
protocol~\eqref{eqn:wit_ave}--\eqref{eqn:bargaininge_ave}, we have 
with probability~1:
\begin{itemize}
\item[(a)]  
The sequence $\{\sum_{i=1}^n \|x^i(t+1)-z\|^2\}$ converges for every 
$z$ in the relative interior of $C(\vmean)$.
\item[(b)] The errors $\bar e^i(t)$ in~\eqref{eqn:error1} are such that
$\sum_{t=0}^\infty\sum_{j=1}^n \|\bar e^i(t)\|^2<\infty$. In particular,
$\lim_{t\to\infty} \|\bar e^i(t)\|=0$ for all $i\in N$.
\end{itemize}
\end{lemma}

\begin{proof}
Let ${\rm rint}Y$ denote the relative interior of a set $Y$. Let $z\in
{\rm rint}C(\vmean)$ be arbitrary and fixed for the rest of the proof.
By Lemma~\ref{lemma:relint}, there exists $t_z$ large enough such that 
$z\in{\rm rint}C(\bar v(t))$ for all $t\ge t_z$ with probability~1. 
Since $C(\bar v(t))=\cap_{i=1}^n X_i(t)$, it follows that $z\in \bar X_i(t)$ 
for all $i\in N$ and $t\ge t_z$ with probability~1.

From $x^i(t+1)=\bar w^i(t) + \bar e^i(t)$ (cf.~\eqref{eqn:bargaininge_ave}),
the definition of $\bar e^i(t)$ in~\eqref{eqn:error2}, 
and the projection property given in Eq.~\eqref{eqn:projection}, we have with 
probability~1, for $i\in N$ and all $t\ge t_z$,
\[\|x^i(t+1) - z\|^2 \leq \|\bar w^i(t) - z\|^2 
-\|\bar e^i(t)\|^2.\]
Summing these relations over $i\in N$, and using 
$\bar w^i(t)=\sum_{j=1}^n a_{ij} x^i(t)$, the convexity of the norm 
and the stochasticity of $A(t)$, we obtain
\[\sum_{i=1}^n\|x^i(t+1) - z\|^2
\le \sum_{i=1}^n\sum_{j=1}^n a_{ij}(t)\|x^j(t)-z\|^2 
- \sum_{j=1}^n\|\bar e^i(t)\|^2.\]
Exchanging the order of summation in 
$\sum_{i=1}^n\sum_{j=1}^n a_{ij}(t)\|x^j(t)-z\|^2 $ and using the doubly 
stochasticity of $A(t)$, we have for all $t\ge t_z$ with probability~1,
\[ \sum_{i=1}^n\|x^i(t+1) - z\|^2
\le \sum_{j=1}^n \|x^j(t)-z\|^2 - \sum_{j=1}^n \|\bar e^i(t)\|^2.\]
Applying the super-martingale converge theorem 
(Theorem~\ref{thm:supermartingale}) (with an index shift),
we see that the sequence $\{\sum_{i=1}^n\|x^i(t+1) - z\|^2\}$ is convergent
and $\sum_{t=t_z}^\infty\sum_{j=1}^n \|\bar e^i(t)\|^2<\infty$ with 
probability~1. Hence, the result in part (b) follows.
\end{proof}

We observe that Lemma~\ref{lemma:xityt} applies to 
protocol~\eqref{eqn:wit_ave}--\eqref{eqn:bargaininge_ave} in view of the
analogy of the description of the protocol 
in~\eqref{eqn:wit}--\eqref{eqn:bargaininge} 
and the protocol in~\eqref{eqn:wit_ave}--\eqref{eqn:bargaininge_ave}.
We will re-state this lemma for an easier reference, but without the proof
since it is almost the same as that of Lemma~\ref{lemma:xityt} (the
proof in essence depends mainly on the matrices $A(k)$).

\begin{lemma}\label{lemma:xityt_ave}
Let Assumptions~\ref{assum:weights} and~\ref{assum:graph} hold. Suppose that 
for the bargaining protocol~\eqref{eqn:wit_ave}--\eqref{eqn:bargaininge_ave} 
we have with probability~1,
\[\lim_{t\to\infty} \|\bar e^i(t)\|=0\qquad\hbox{for all }i\in N.\]
Then, for every player $i\in N$ we have with probability~1,
\[\lim_{t\to\infty}\|x^i(t)-y(t)\|=0,\qquad
\lim_{t\to\infty}\|\bar w^i(t)-y(t)\|=0,\]
where $y(t)=\frac{1}{n}\sum_{j=1}^n x^j(t)$.
\end{lemma}

We are now ready to show the convergence of the bargaining protocol.
We show this in the following theorem by combining the basic properties 
of the protocol established in Lemma~\ref{lemma:basic_ave} and 
using Lemma~\ref{lemma:xityt_ave}.

\begin{theorem}\label{thm:conv_ave}
Let Assumptions~\ref{assum:weights}--\ref{assum:cores} hold.
Then, the bargaining protocol~\eqref{eqn:wit_ave}--\eqref{eqn:bargaininge_ave} 
converges to a random allocation that lies in the core $C(\vmean) $ of 
the average game with probability~1, i.e., there is a random vector 
$\tilde z\in C(\vmean)$  such that with probability~1,
$$\lim_{t\rightarrow \infty} \|x^i(t)-\tilde z\|=0 \qquad 
\mbox{for all $i\in N$}.$$
\end{theorem}

\begin{proof} 
By Lemma~\ref{lemma:basic_ave}, with probability~1, the sequence 
$\{\sum_{i=1}^n \|x^i(t)-z\|^2\}$ is convergent for 
every $z$ in the relative interior of $C(\vmean)$ and $\|e^i(t)\|\to0$ 
for each player $i\in N$. By Lemma~\ref{lemma:xityt_ave}, with probability~1 
we have 
\begin{equation}\label{eqn:xtyt_ave}
\lim_{t\to\infty} \|x^i(t)-y(t)\|=0\qquad \hbox{for every $i\in N$}.
\end{equation} 
Hence, with probability~1,
\begin{equation}\label{eqn:yt}
\hbox{$\{\|y(t)-z\|\}$ is convergent for 
every $z\in {\rm rint}C(\vmean)$}.\end{equation}

We next show that $\{y(t)\}$ has accumulation points in the core $C(\vmean)$ 
with probability~1. Since $x^i(t+1) \in \bar X_i(t)$, it follows
\begin{equation}\label{eqn:ytx}
\sum_{i=1}^n \dist^2\left(y(t+1),\bar X_i(t)\right) 
\leq \sum_{i=1}^n \|y(t+1) - x^i(t+1)\|^2.
\end{equation} 
Relations~\eqref{eqn:ytx} and~\eqref{eqn:xtyt_ave} 
imply 
\begin{align}\label{eqn:lim}
\lim_{t\rightarrow \infty} 
\sum_{i=1}^n \dist^2\left(y(t+1),\bar X_i(t)\right)=0 \quad
\hbox{with probability~1}.
\end{align}
Let $i\in N$ be arbitrary but fixed and let $x\in \bar X_i$ be arbitrary. 
Note that we can write
\[\|y(t+1) - x\|\le \|y(t+1) - P_{\bar X_i(t)} [y(t+1)]\| +
\|P_{\bar X_i(t)} [y(t+1)]-x\|.\]
By taking the minimum over $x\in\bar X_i$ and 
using $\|y(t+1)-P_{\bar X_i(t)} [y(t+1)]\|=\dist(y(t+1),\bar X_i(t))$,
we obtain for all $i\in N$ and $t\ge0$,
\[\dist(y(t+1),\bar X_i)\le \dist(y(t+1),\bar X_i(t)) +
\dist\left(P_{\bar X_i(t)} [y(t+1)], \bar X_i\right).\]
The bounding sets $\bar X_i$ are nonempty by Assumption~\ref{assum:cores}
and the fact $C(\vmean)\subset \bar X_i$ for all $i$,
while $\bar X_i(t)$ are assumed nonempty (see the discussion 
after relation~\eqref{eqn:bs_ave}).
Furthermore, since $\bar X_i$ and $\bar X_i(t)$ are polyhedral sets,
by Lemma~\ref{lemma:twopolyhsets} it follows that for a scalar $L_i>0$ and all
$t\ge0$,
\[\dist\left(P_{\bar X_i(t)} [y(t+1)], \bar X_i\right)
\le L_i \|\bar v(t) - \vmean\|.\]
Therefore, for all $i\in N$ and all $t\ge0$,
\[\dist(y(t+1),\bar X_i)\le \dist(y(t+1),\bar X_i(t)) 
+ L_i \|\bar v(t) - \vmean\|.\]
By letting $t\to\infty$ and using relations~\eqref{eqn:lim} 
and $\bar v(t)\to\vmean$ (Assumption~\ref{assum:ergodic}), we see
that for all $i\in N$ with probability~1,
\begin{align}\label{eqn:alli}
\lim_{t\to\infty}\,\dist(y(t+1),\bar X_i)=0.\end{align}

In view of relation~\eqref{eqn:yt},
the sequence $\{y(t)\}$ is bounded with probability~1,
so it has accumulation points with probability~1. By relation~\eqref{eqn:alli},
all accumulation points of $\{y(t+1)\}$ lie in the set $\bar X_i$ for every 
$i\in N$ with probability~1. Therefore, the accumulation points of 
$\{y(t)\}$ must lie in the intersection $\cap_{i\in N} \bar X_i$ with 
probability~1. Since $\cap_{i\in N} \bar X_i=C(\vmean)$, we conclude that all 
accumulation points of $\{y(t)\}$ lie in the core $C(\vmean)$ 
with probability~1. Furthermore, according to relation~\eqref{eqn:yt} we have
that, for any point $z\in{\rm rint}C(\vmean)$,  the accumulation points of 
the sequences $\{y(t)\}$ are at the same (random) distance from $z$ 
with probability~1. Since the accumulation points are in the set $C(\vmean)$, 
it follows that $\{y(k)\}$ is convergent with probability~1 and its limit point
is in the core $C(\vmean)$ with probability~1. Now, since $\|x^i(t)-y(t)\|\to0$
with probability~1 for all $i$ (see~\eqref{eqn:xtyt_ave}), the sequences 
$\{x^i(t)\}$, $i\in N$, have the same limit point as the sequence $\{y(t)\}$. 
Thus, the sequences $\{x^i(t)\}$, $i\in N$, converge to a common (random) point
in $C(\vmean)$ with probability~1.
\end{proof}

Theorem~\ref{thm:conv_ave} shows the convergence of the allocations generated 
by the bargaining protocol in~\eqref{eqn:wit_ave}--\eqref{eqn:bargaininge_ave}.
The convergence relies on the properties of the matrices and the connectivity 
of the players' neighbor graphs, as reflected in 
Assumptions~\ref{assum:weights} and~\ref{assum:graph}. It also critically
depends on the fact that the core of the average game $C(\vmean)$ has
dimension $n-1$ and that all bounding sets $\bar X_i(t)$ and, hence
the cores $C(\bar v(t))$, lie in the same hyperplane with probability~1,
the hyperplane defined through the constant value of the grand coalition.

\section{Numerical Illustrations}\label{sec:numerical}
In this section, we report some numerical simulations.
We consider coalitional TU games with 3 players, so the number of possible 
nonempty coalitions is $m=7$. We consider two different scenarios as shown 
respectively in rows I and II in Table~\ref{Tabb2}. The columns
of Table~\ref{Tabb2} enumerate the coalitions. In each scenario,
the characteristic functions $v_S(t)$ are generated independently with 
identical uniform distribution over an interval. Specifically, we suppose that 
the values of single player coalitions $\{1\}$ and $\{2\}$ are uncertain
within the given interval. All the other coalitions' values are fixed and
equal to zero except for the grand coalition, which has value~10 in both 
scenarios~I and~II.

\begin{table}[b]
	\centering
		\begin{tabular}{|c|c|c|c|c|c|c|c|}
\hline		  
& $v_{\{1\}}$	& $v_{\{2\}}$ & $v_{\{3\}}$ & $v_{\{1,2\}}$ & $v_{\{1,3\}}$ & 
$v_{\{2,3\}}$ & $v_{\{1,2,3\}}$ \\\hline
I & 	 $[4,7]$ &  $[0,3]$ & $0$ & $0$ & $0$ & $0$ & $10$ \\\hline
II & 	 $[4,9]$ &  $[0,5]$ & $0$ & $0$ & $0$ & $0$ & $10$ \\\hline
		\end{tabular}
\caption{Coalitions' values for the two simulations scenarios.}\label{Tabb2}
\end{table}

The two scenarios differ in that the core $C(v^{\max})$ of the robust game
is nonempty in scenario I and empty in scenario II.
For scenario I, we simulate the convergence behavior of
the bargaining protocol~\eqref{eqn:wit}--\eqref{eqn:bargaininge} for the 
robust game, while for scenario II, we simulate bargaining 
protocol~\eqref{eqn:wit_ave}--\eqref{eqn:bargaininge_ave} for the average game.

\begin{figure}[b!]
\centering
\subfigure[]{
\includegraphics[scale=.6]{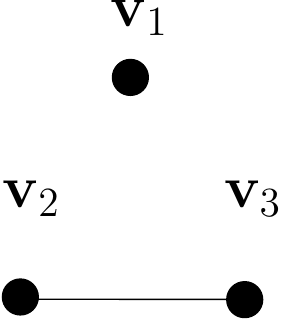}
\label{fig:subfig11}}
\hspace{1.5cm}
\subfigure[]{
\includegraphics[scale=.6]{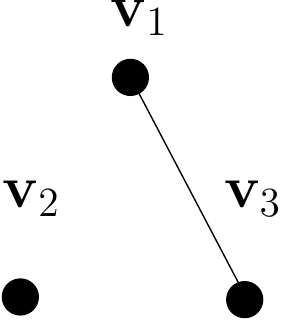}
\label{fig:subfig22}}
\hspace{1.5cm}
\subfigure[]{
\includegraphics[scale=.6]{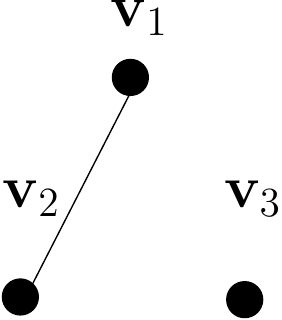}
\label{fig:subfig33}
}
\label{fig:subfigureExample1}
\caption{\small Topology of players' neighbor-graph 
at three distinct times $t=0$, 1 and 2.}
\end{figure}

For each scenario, we run 50 different Monte Carlo trajectories each one having
$100$ iterations. The number of iterations is chosen long enough to show 
the convergence of the protocols. All plots include the sampled average 
and sampled variance for the 50 different trajectories that were simulated.
Each trajectory in each scenario is generated by starting with the same initial
allocations, which are given by
$x^1(0)=[10 \; 0 \; 0]'$, $x^2(0)=[0 \; 10 \; 0]'$, and 
$x^3(0)=[0 \; 0 \; 10]'$. The sampled average is computed for each time 
$t=1,\ldots,100$, by fixing the time $t$ and computing the average value 
of the 50 trajectory sample values for that time. The sampled variance is 
computed as the variance of the samples with respect to their sampled average.

Regarding the players' neighbor-graphs, we assume that the graphs are 
deterministic but time-varying. The graphs for the times $t=0,1,2$ are 
as follows: player 2 and 3 connected at time $t=0$ 
(see Figure~\ref{fig:subfig11}),
then player 3 and 1 connected at time $t=1$ (Figure~\ref{fig:subfig22}), 
and finally  player 1 and 2 connected at time $t=3$ 
(Figure~\ref{fig:subfig33}). These graphs are then repeated consecutively
in the same order. In this way, the players' neighbor-graph is connected every 
3 time units (Assumption~\ref{assum:graph} is satisfied with $Q=2$).

The matrices that we associate with these three graphs, are respectively 
given by:
$$A(0)=\left[\begin{array}{ccc} 1 & 0 & 0\\0 & \frac{1}{2} & \frac{1}{2}\\
0 & \frac{1}{2} & \frac{1}{2}\end{array}\right],\qquad
A(1)=\left[\begin{array}{ccc} \frac{1}{2}  & 0 & \frac{1}{2} \\
0 & 1 & 0\\\frac{1}{2}  & 0 & \frac{1}{2}\end{array}\right],\qquad
A(2)=\left[\begin{array}{ccc} \frac{1}{2} & \frac{1}{2} & 0\\
\frac{1}{2} & \frac{1}{2} & 0\\0 & 0 & 1\end{array}\right].$$
These matrices are also repeated in the same order for the rest of the time.
Thus, at any time $t$, the matrix $A(t)$ is doubly stochastic, with positive
diagonal, and every positive entry bounded below by $\frac{1}{2}$, so 
Assumption~\ref{assum:weights} is satisfied with $\alpha=\frac{1}{2}$.
All simulations are carried out with MATLAB on an Intel(R) Core(TM)2 Duo, 
CPU P8400 at 2.27 GHz and a 3GB of RAM. The run time of each simulation
is around 90 seconds.
\subsection{Simulation Scenario I}
In this scenario, the coalitions' values are generated as given in row~I of 
Table~\ref{Tabb2}. In particular, at each time~$t$, the value $v_{\{1\}}(t)$ is
chosen randomly in the interval $[4,7]$ with uniform probability independently 
of the other times. Similarly, the values $v_{\{2\}}(t)$ are generated in the 
interval $[0,3]$. The grand coalition value is fixed to 10 at all times, and 
the other coalition values are 0. With this data, we consider the allocations 
as generated by players $i=1,2,3$ according to the bargaining protocol
in~\eqref{eqn:wit}--\eqref{eqn:bargaininge} 
for the robust game as reported in Subsection~\ref{sec:simI_robust}.
We then consider the bargaining protocol 
in~\eqref{eqn:wit_ave}--\eqref{eqn:bargaininge_ave} for the average game 
in Subsection~\ref{sec:simI_aver}.

\subsubsection{Robust game}\label{sec:simI_robust}
For this specific example, the characteristic function $v^{\max}$ for 
the robust game is obtained by considering the highest possible coalition 
values (see the data in row~I of Table~\ref{Tabb2}), which results in 
$v^{\max}=[7\,3\,0\,0\,0\,0\,10]'$. The resulting core of the robust game is 
given by
\begin{eqnarray*}C(v^{\max})&=&\{x\in \mathbb R^3:\, x_1 \geq 7,\, x_2\geq
3,\,x_3\geq 0,\, x_1+x_2\geq 0,\,x_1+x_3\geq 0,\\&& x_2+x_3\geq
0,\,x_1+x_2+x_3=10\}.\end{eqnarray*}
We note that this core contains a single point, namely  $[7\,3\,0]'$.

To ensure that $v(t)=v^{\rm max}$ infinitely often, 
as required by Theorem~\ref{thm1} for the convergence of the protocol,
we adopt the following randomization mechanism. At each time 
$t=1,\ldots,100$, we flip a coin and if the outcome is ``head'' 
(probability 1/2), the coalitions' values 
$v_{\{1\}}(t)$ and $v_{\{2\}}(t)$ are extracted from the intervals $[4,7]$ and 
$[0,3]$, respectively, with uniform probability independently of the other 
times. If  the outcome of the coin flip is ``tail'', then we assume that 
the robust game realizes and take $v(t)=v^{\rm max}$. 

We next present the results obtained by the Monte Carlo runs for the 
bargaining protocol in~\eqref{eqn:wit}--\eqref{eqn:bargaininge}.
An illustration of a typical run with the allocations generated in periods 
$t=0,1,2,3$ is shown below:
\begin{eqnarray*}
\begin{array}{l} 
v(0)=[6.8 \; 2.7 \ldots 10]'\\
v(1)=[7 \; 3 \ldots 10]'\\
v(2)=[4.4 \; 1.1 \ldots 10]'\\
v(3)=[7 \; 3 \ldots 10]'\end{array}
\,\begin{array}{l}
x^1(0)=[10 \; 0 \; 0]'\\
x^1(1)=[10 \; 0 \; 0]'\\
x^1(2)=[5 \; 2.5 \; 2.5]'\\
x^1(3)=[7 \; 1.5 \; 1.5]'\end{array}
\, \begin{array}{l}
x^2(0)=[0 \; 10 \; 0]'\\
x^2(1)=[0 \; 5 \; 5]'\\
x^2(2)=[0 \; 5 \; 5]'\\
x^2(3)=[2.5 \; 3.75 \; 3.75]'\end{array}
\, \begin{array}{l}
x^3(0)=[0 \; 0 \; 10]'\\
x^3(1)=[0 \; 5 \; 5]'\\
x^3(2)=[5 \; 2.5 \; 2.5]'\\
x^3(3)=[5 \; 2.5 \; 2.5]'.\end{array}\end{eqnarray*} 
Recall that the initial allocations of the players are 
$x^1(0)=[10 \; 0 \; 0]'$, $x^2(0)=[0 \; 10 \; 0]'$, and 
$x^3(0)=[0 \; 0 \; 10]'$. At time $t=1$, bargaining involves player~2 and~3 
who update the allocations respectively as  $x^2(1)=[0 \; 5 \; 5]'$ and 
$x^3(1)=[0 \; 5 \; 5]'$. These allocations are feasible for their bounding sets
so the projections on these sets are not performed. At time $t=2$, 
the bargaining involves player~1 and~3 who update their allocations, 
respectively, as  $x^1(2)=[5 \; 2.5 \; 2.5]'$ and $x^3(2)=[5 \; 2.5 \; 2.5]'$. 
Again, these allocations are feasible for their bounding sets and 
the projections are not performed. Finally, at time $t=3$, the bargaining 
involves player~1 and~2 who update their allocations resulting in  
$x^1(3)=[7 \; 1.5 \; 1.5]'$ and $x^2(3)=[2.5 \; 3.75 \; 3.75]'$. 
Notice that $x^1(3)$ is obtained after player~1 projects onto his bounding set.

\begin{figure}
\centering
\includegraphics[width=15cm]{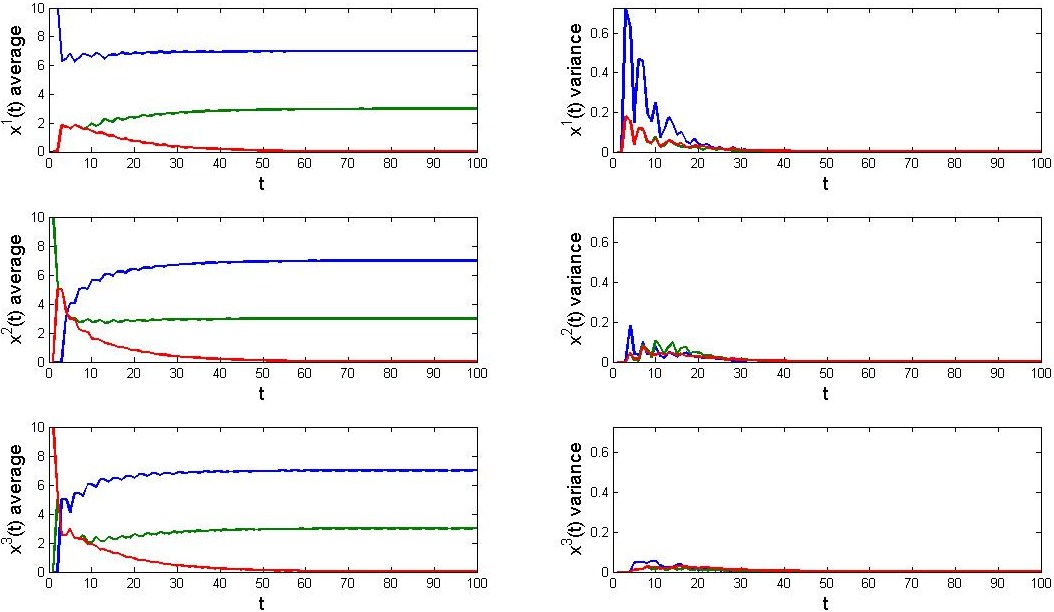}
\caption{Plots of the sampled average (left) and variance (right) of 
players' allocations $x^i(t)$, $i=1,2,3$ generated by bargaining 
protocol~\eqref{eqn:wit}--\eqref{eqn:bargaininge} for 
the robust game associated with the data in row I of Table~\ref{Tabb2}. 
Sampled averages of the allocations
$x^i(t)$ converge to the same point $\tilde x=[7 \, 3 \, 0]'\in C(v^{\max})$,
while sampled variances go rapidly to zero.}  \label{p:Inst1}
\end{figure}

\begin{figure}
\centering
\includegraphics[width=15cm]{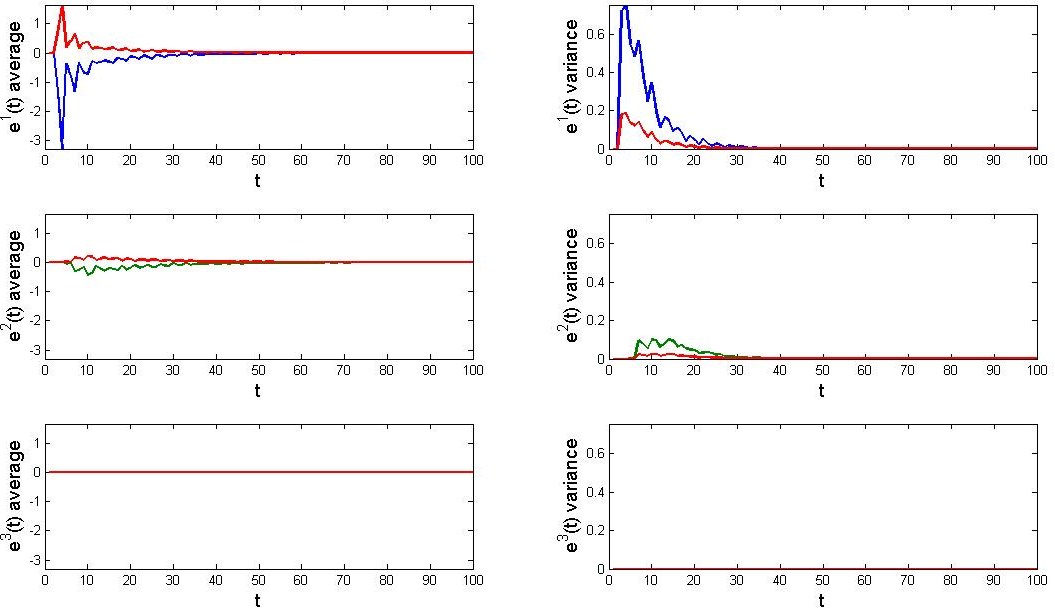}
\caption{Plots of the sampled average (left) and sampled variance (right) of 
the errors $e^i(t)$, $i=1,2,3$ for the bargaining 
protocol~\eqref{eqn:wit}--\eqref{eqn:bargaininge} and the robust
game associated with the data in row I of Table~\ref{Tabb2}. 
Sampled averages and the variances of the errors $e^i(t)$ converge to zero.}
\label{p:Inst3}
\end{figure}

In Figures~\ref{p:Inst1} and~\ref{p:Inst3}, 
we report our simulation results for the average of the sample 
trajectories obtained by Monte Carlo runs. 
Figure~\ref{p:Inst1} shows the sampled average and variance of the allocations 
$x^i(t)$, $i=1,2,3$ per iteration~$t$. In accordance with the convergence 
result of Theorem~\ref{thm1}, the sampled averages of the players' allocations 
$x^i(t)$ converge to the same point, namely $x=[7 \, 3 \, 0]'$ which is in 
the core of the robust game $C(v^{\rm max})$. Figure~\ref{p:Inst3} shows that 
the sample average and sampled variance of the errors $e^i(t)$ converge to 0, 
as expected in view of Lemma~\ref{lemma:basic}(b). 

\subsubsection{Average game}\label{sec:simI_aver}

In this numerical example, for scenario I data as given in row~I of 
Table~\ref{Tabb2}, we consider the average TU game and its corresponding 
bargaining protocol~\eqref{eqn:wit_ave}--\eqref{eqn:bargaininge_ave}.
We at first provide a sample example and report the outcomes for the
first three steps of the algorithm, as reported below:
\begin{eqnarray*} 
\begin{array}{l} \bar v(0)=[5.5 \; 1.6 \ldots 10]'\\
\bar v(1)=[4.8 \; 1.8 \ldots 10]'\\
\bar v(2)=[4.6 \; 1.8 \ldots 10]'\\
\bar v(3)=[4.7 \; 1.8 \ldots 10]'\end{array}
\,\begin{array}{l}
x^1(0)=[10 \; 0 \; 0]'\\
x^1(1)=[10 \; 0 \; 0]'\\
x^1(2)=[5 \; 2.5 \; 2.5]'\\
x^1(3)=[4.7 \; 2.6 \; 2.7]'\end{array}
\, \begin{array}{l}
x^2(0)=[0 \; 10 \; 0]'\\
x^2(1)=[0 \; 5 \; 5]'\\
x^2(2)=[0 \; 5 \; 5]'\\
x^2(3)=[2.5 \; 3.75 \; 3.75]'\end{array}
\, \begin{array}{l}
x^3(0)=[0 \; 0 \; 10]'\\
x^3(1)=[0 \; 5 \; 5]'\\
x^3(2)=[5 \; 2.5 \; 2.5]'\\
x^3(3)=[5 \; 2.5 \; 2.5]'.\end{array}\end{eqnarray*}
Recalling that the initial players' allocations are $x^1(0)=[10 \; 0 \; 0]'$, 
$x^2(0)=[0 \; 10 \; 0]'$, and $x^3(0)=[0 \; 0 \; 10]'$,
we note that at time $t=1$, bargaining involves player~2 and~3 who update their
allocations to $x^2(1)=[0 \; 5 \; 5]'$ and $x^3(1)=[0 \; 5 \; 5]'$, 
respectively. At time $t=2$, the bargaining involves player~1 and~3 who update 
their allocations, respectively, to $x^1(2)=[5 \; 2.5 \; 2.5]'$ and 
$x^3(2)=[5 \; 2.5 \; 2.5]'$. Notice that the obtained allocations for $t=1,2$ 
are feasible for the bounding sets and, therefore, no projection is performed.
Finally, at time $t=3$, the bargaining involves player~1 and~2 who update 
their allocations to $x^1(3)=[4.7 \; 2.6 \; 2.7]'$ and 
$x^2(3)=[2.5 \; 3.75 \; 3.75]'$, respectively. Here $x^1(3)$ results from 
projecting onto the bounding set of player~1.

For the average game associated with the data in row I of 
Table~\ref{Tabb2}, we have $v^{\rm mean}=[5.5\ 1.5\ 0\ 0\ 0\ 0\ 10]'$ and 
the core $C(v^{\rm mean})$ given by
\begin{eqnarray*}C(v^{\rm mean})=\{x\in \mathbb R^3:\, x_1 \geq 5.5,\, x_2\geq
1.5,\,x_3\geq 0,\, x_1+x_2\geq 0,\,x_1+x_3\geq 0,\\x_2+x_3\geq
0,\,x_1+x_2+x_3=10\}.\end{eqnarray*}

\begin{figure} 
\centering
\includegraphics[width=15cm]{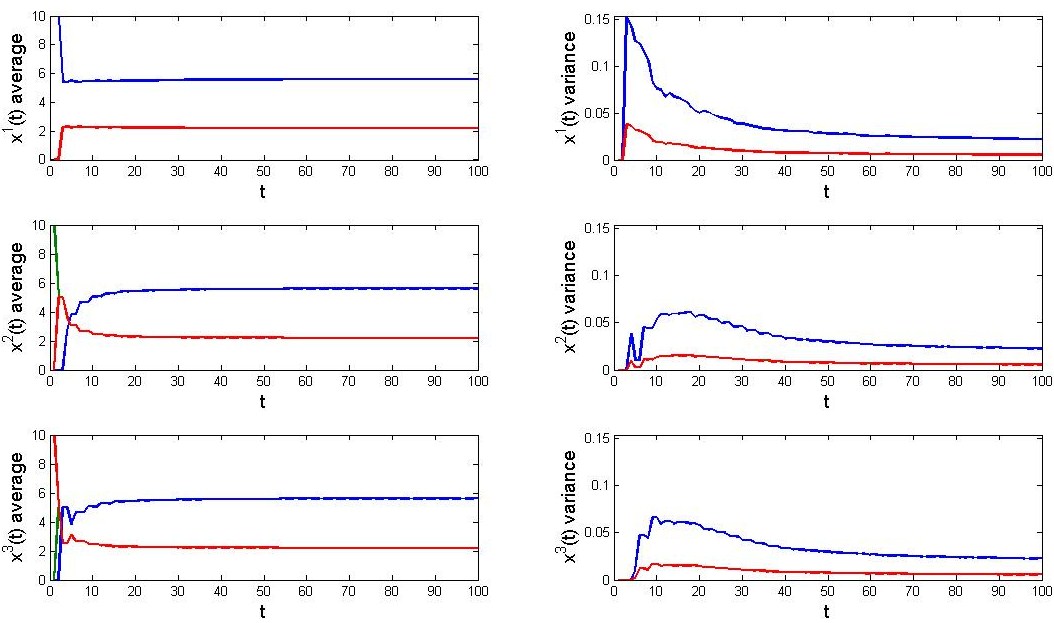}
\caption{Plots of the sampled averages (left) and variances (right) of 
players' allocations $x^i(t)$, $i=1,2,3,$ obtained by bargaining 
protocol~\eqref{eqn:wit_ave}--\eqref{eqn:bargaininge_ave} for the average 
game associated with the data in row I of Table~\ref{Tabb2}. Sampled 
averages converge to point $\tilde x=[5.6 \ 2.2 \ 2.2]' \in C(v^{\rm mean})$,
which is the average of the limit points of the 50-sample trajectories.
}\label{p:Ave1}
\end{figure}

\begin{figure} 
\centering
\includegraphics[width=15cm]{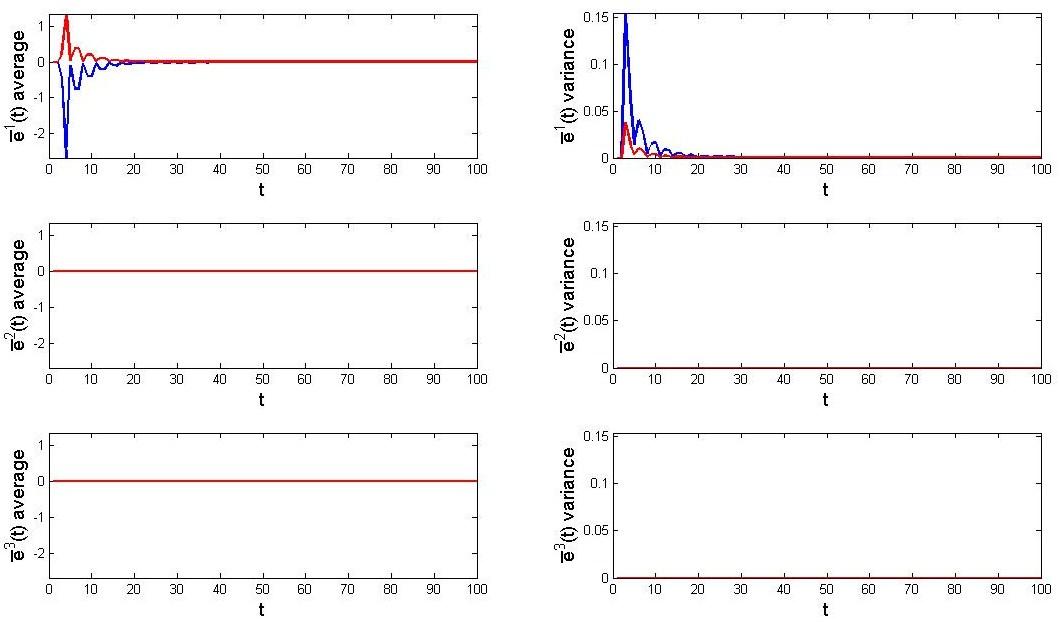}
\caption{Plots of the sampled average (left) and variance (right) of 
the errors $\bar e^i(t)$, $i=1,2,3,$ for bargaining 
protocol~\eqref{eqn:wit_ave}--\eqref{eqn:bargaininge_ave} and the average game 
associated with the data in row~I of Table~\ref{Tabb2}. 
The sampled average and variance of the errors $\bar e^i(t)$ goes to~0.}
\label{p:Ave3}
\end{figure}

In Figures~\ref{p:Ave1} and~\ref{p:Ave3}, we depict our simulation results
generated by bargaining 
protocol~\eqref{eqn:wit_ave}--\eqref{eqn:bargaininge_ave}
for the average game. Figure~\ref{p:Ave1} shows that the sampled average of 
the allocations $x^i(t)$, $i=1,2,3$ converge to a common point
$\tilde x=[5.6 \ 2.2 \ 2.2]'$ which belongs to the core $C(\vmean)$
of the average game, as guaranteed by Theorem \ref{thm:conv_ave}. 
The sampled variance does not converge to zero as the common limit point of 
the allocations $x^i(t)$ can be different for different runs. 
Figure~\ref{p:Ave3} demonstrates that the sampled average and variance of 
the errors  $\bar e^i(t)$, $i=1,2,3,$ converge to zero, as predicted by
Lemma~\ref{lemma:basic_ave}(b).   

\subsection{Simulation Scenario II}
\begin{figure} [h!]
\centering
\includegraphics[width=15cm]{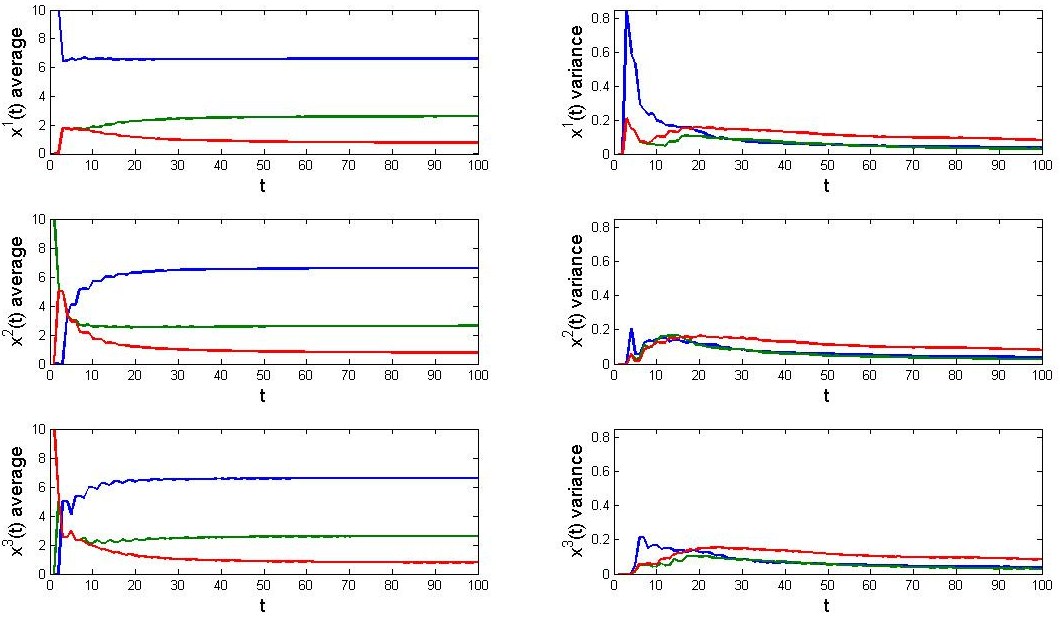}
\caption{Plots of the sampled averages (left) and variances (right) of 
the players' allocations $x^i(t)$, $i=1,2,3$ obtained by the bargaining 
protocol~\eqref{eqn:wit_ave}--\eqref{eqn:bargaininge_ave}
for the average game associated with the data in row II of Table~\ref{Tabb2}.
Sampled averages converge to point 
$\tilde x=[6.6 \ 2.6 \ 0.8]' \in C(v^{\rm mean})$, which is the average
of the limit points of the 50-sample trajectories.}
\label{p:AveEmptyCore28082010-1}
\end{figure}

\begin{figure} [h!]
\centering
\includegraphics[width=15cm]{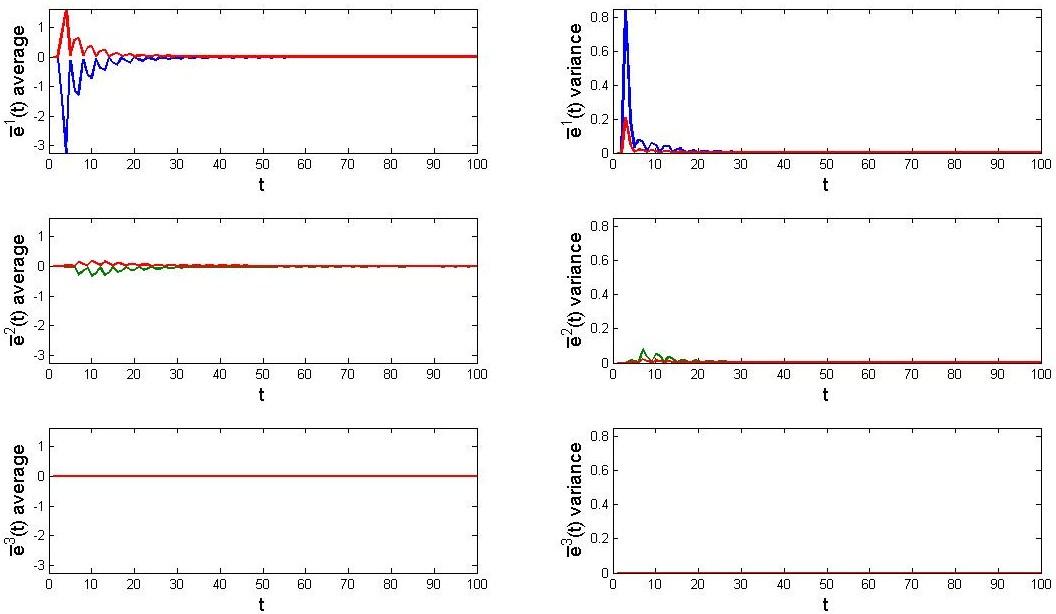}
\caption{Plots of the sampled average and variance of the errors 
$\bar e^i(t)$, $i=1,2,3$ for the bargaining 
protocol~\eqref{eqn:wit_ave}--\eqref{eqn:bargaininge_ave} 
and the average game associated with the data in row~II of Table~\ref{Tabb2}. 
The sample average and variance of the errors  $\bar e^i(t)$ converge to zero.}
\label{p:AveEmptyCore28082010-3}
\end{figure}

Here, we report the simulation results obtained 
by the bargaining protocol~\eqref{eqn:wit_ave}--\eqref{eqn:bargaininge_ave} for
the average game corresponding to the data in row~II of Table~\ref{Tabb2}.
In this case the core of the robust game is empty, so we do not consider the 
robust game. The average game has characteristic function 
$v^{\rm mean}=[6.5\,2.5\,0\,0\,0\,0\,10]'$ and its core is 
\begin{eqnarray*}C(v^{\rm mean})=\{x\in \mathbb R^3:\, x_1 \geq 6.5,\, x_2\geq
2.5,\,x_3\geq 0,\, x_1+x_2\geq 0,\,x_1+x_3\geq 0,\\x_2+x_3\geq
0,\,x_1+x_2+x_3=10\}.\end{eqnarray*}
Figures~\ref{p:AveEmptyCore28082010-1} and \ref{p:AveEmptyCore28082010-3}
show the results for the average game obtained in our simulations.
In Figure~\ref{p:AveEmptyCore28082010-1}, we report the sampled averages 
of the players' allocations $x^i(t)$, $i=1,2,3$, obtained by bargaining 
protocol~\eqref{eqn:wit_ave}--\eqref{eqn:bargaininge_ave}.
In accordance with Theorem~\ref{thm:conv_ave}, the players' allocations 
converge to an allocation that lies in the core of the average game 
$C(v^{\rm mean})$, precisely to the point 
$\tilde x=[6.6 \ 2.6 \ 0.8]' \in C(v^{\rm mean})$. Here, again,
the sampled variance of the allocations does not converge to zero as 
their common limit point is different for different runs.
Figure~\ref{p:AveEmptyCore28082010-3} shows that the sampled average and 
variance of the errors  $\bar e^i(t)$, $i=1,2,3$, converge to zero.

\section{Conclusions}\label{conclusions}
This article deals with dynamics and robustness within the framework
of coalitional TU games. With respect to a sequence of TU games,
each with a random characteristic function, 
the novelty of the work lies in the design of a decentralized
allocation process defined over a communication graph of players. 
The allocation process captures the main features
of bargaining in a realistic scenario. The proposed bargaining scheme is 
proven to converge, with probability~1, in either a robust game setting or 
an average game setting, under mild assumptions on the communication 
topology among the players and  
the stochastic properties of the random characteristic function.
The key properties that distinguish this work from the existing work
on dynamic games are: (1) the introduction of a time-varying communication 
graph, termed players' neighbor-graph, over which the bargaining protocol
takes place; and (2) the distributed bargaining protocol for 
players' allocations updates subject to local information exchange
with neighboring players.
 
\bibliographystyle{plain}        
\bibliography{biblio}           

\begin{thebibliography}{10}

\bibitem{AS02}
T.~Arnold and U.~Schwalbe.
\newblock Dynamic coalition formation and the core.
\newblock {\em Journal of Economic Behavior and Organization}, 49:363--380,
  2002.

\bibitem{BBP10}
D.~Bauso, F.~Blanchini, and R.~Pesenti.
\newblock Optimization of long-run average-flow cost in networks with
  time-varying unknown demand.
\newblock {\em IEEE Transactions on Automatic Control}, 55(1):20--31, 2010.

\bibitem{BR10}
D.~Bauso and P.~V. Reddy.
\newblock Robust allocation rules in dynamical cooperative {TU} games.
\newblock {\em Proc. of the 49th Conference on Decision and Control}, 2010.

\bibitem{BT09}
D.~Bauso and J.~Timmer.
\newblock Robust dynamic cooperative games.
\newblock {\em International Journal of Game Theory}, 38(1):23--36, 2009.

\bibitem{C98}
J.C. Cesco.
\newblock A convergent transfer scheme to the core of a {TU}-game.
\newblock {\em Revista de Matem\'aticas Aplicadas}, 19(1--2):23--35, 1998.

\bibitem{Drechsel10}
J.~Drechsel.
\newblock {\em Cooperative Lot Sizing Games in Supply Chains.}
\newblock Springer-Verlag. Series: Lecture Notes in Economics and Mathematical
  Systems, Berlin, Germany, 1 edition, 2010.

\bibitem{Facchinei2003}
F.~Facchinei and J-S. Pang.
\newblock {\em Finite-dimensional variational inequalities and complementarity
  problems}, volume I-II.
\newblock Springer-Verlag, New York, 2003.

\bibitem{FP00}
J.A. Filar and L.A. Petrosjan.
\newblock Dynamic cooperative games.
\newblock {\em International Game Theory Review}, 2(1):47--65, 2000.

\bibitem{G77}
D.~Granot.
\newblock Cooperative games in stochastic characteristic function form.
\newblock {\em Management Sci.}, 23:621--630, 1977.

\bibitem{HDS00}
B.C. Hartman, M.~Dror, and M.~Shaked.
\newblock Cores of inventory centralization games.
\newblock {\em Games and Economic Behavior}, 31:26--49, 2000.

\bibitem{H75}
A.~Haurie.
\newblock On some properties of the characteristic function and the core of a
  multistage game of coalitions.
\newblock {\em IEEE Transactions on Automatic Control}, 20(2):238--241, 1975.

\bibitem{Hoffman1952}
A.J. Hoffman.
\newblock On approximate solutions of systems of linear inequalities.
\newblock {\em Journal of Research of the National Bureau of Standards},
  49:263--265, 1952.

\bibitem{L02}
E.~Lehrer.
\newblock Allocation processes in cooperative games.
\newblock {\em International Journal of Game Theory}, 31:341--351, 2002.

\bibitem{NOOT08}
A.~Nedi\'c, A.~Olshevsky, A.~Ozdaglar, and J.N. Tsitsiklis.
\newblock Distributed subgradient methods and quantization effects.
\newblock {\em Proc. of the 47th CDC Conference}, 54(11):4177--4184, 2008.

\bibitem{NOOT09}
A.~Nedi\'c, A.~Olshevsky, A.~Ozdaglar, and J.N. Tsitsiklis.
\newblock On distributed averaging algorithms and quantization effects.
\newblock {\em IEEE Transactions on Automatic Control}, 54(11):2506--2517,
  2009.

\bibitem{Nedic09}
A.~Nedi\'c and A.~Ozdaglar.
\newblock Distributed subgradient methods for multi-agent optimization.
\newblock {\em IEEE Transactions on Automatic Control}, 54(1):48--61, 2009.

\bibitem{NOP10}
A.~Nedi\'c, A.~Ozdaglar, and P.A. Parrilo.
\newblock Constrained consensus and optimization in multi-agent networks.
\newblock {\em IEEE Transactions on Automatic Control}, 55(4):922--938, 2010.

\bibitem{Polyak87}
B.~T. Polyak.
\newblock {\em Introduction to Optimization}.
\newblock Optimization Software, Inc., New York, 1987.

\bibitem{Sundhar09}
S.~Sundhar Ram, A.~Nedi\'c, and V.V. Veeravalli.
\newblock Incremental stochastic subgradient algorithms for convex
  optimization.
\newblock {\em SIAM Journal on Optimization}, 20(2):691--717, 2009.

\bibitem{Sundhar08c}
S.~Sundhar Ram, A.~Nedi\'{c}, and V.V. Veeravalli.
\newblock Distributed stochastic subgradient algorithm for convex optimization.
\newblock {\em Journal of Optimization Theory and Applications},
  147(3):516--545, 2010.

\bibitem{Robbins1971}
H.~Robbins and D.~Siegmund.
\newblock A convergence theorem for nonnegative almost supermartingales and
  some applications.
\newblock In {\em Optimizing Methods in Statistics}, pages 233--257. Academic
  Press, New York, 1971.

\bibitem{Rockafellar1998}
R.~T. Rockafellar and R.~J-B. Wets.
\newblock {\em Variational Analysis}.
\newblock Springer-Verlag, Berlin, Germany, 1998.

\bibitem{SHDHB09}
W.~Saad, Z.~Han, M.~Debbah, A.~Hj{\"o}rungnes, and T.~Ba\c{s}ar.
\newblock Coalitional game theory for communication networks.
\newblock {\em IEEE Signal Processing Magazine, Special Issue on Game Theory},
  26(5):77--97, 2009.

\bibitem{SB99}
J.~Suijs and P.~Borm.
\newblock Stochastic cooperative games: Superadditivity, convexity, and
  certainty equivalents.
\newblock {\em Games and Economic Behavior}, 27(2):331--345, 1999.

\bibitem{TBT03}
J.~Timmer, P.~Borm, and S.~Tijs.
\newblock On three shapley-like solutions for cooperative games with random
  payoffs.
\newblock {\em International Journal of Game Theory}, 32:595--613, 2003.

\bibitem{Tsitsiklis84}
J.N. Tsitsiklis.
\newblock {\em Problems in Decentralized Decision Making and Computation}.
\newblock PhD thesis, Dept. of {E}lectrical {E}ngineering and {C}omputer
  {S}cience, {MIT}, 1984.

\bibitem{VM}
J.~von Neumann and O.~Morgenstern.
\newblock {\em Theory of Games and Economic Behavior}.
\newblock Princeton Univ. Press, 1944.

\end{thebibliography}

\end{document}